\documentclass{amsart}

\usepackage{amsmath, amsthm, amssymb, latexsym,url}

\usepackage{hyperref}

\usepackage{graphicx}
\usepackage{caption}
\usepackage{subcaption}
\usepackage{bbm}

\usepackage{rotating}

\title[Invariant measure for Hurwitz CFs]{Calculations of the invariant measure for Hurwitz Continued Fractions}

\author[G. Hiary]{Ghaith Hiary}
\address{
Department of Mathematics\\
The Ohio State University\\
Columbus, OH 432100}
\email{hiary.1@osu.edu}

\author[J. Vandehey]{Joseph Vandehey}
\address{
Department of Mathematics\\
The Ohio State University\\
Columbus, OH 432100}
\email{vandehey.1@osu.edu}

\date{\today}



\newtheorem{thm}{Theorem}[section]

\newtheorem{conj}[thm]{Conjecture}

\newcommand{\ii}{\mathbbm{i}}

\begin{document} 

\maketitle

\begin{abstract}
We study the density of the invariant measure of the Hurwitz complex continued fraction from a computational perspective. It is known that this density is piece-wise real-analytic and so we provide a method for calculating the Taylor coefficients around certain points and also the results of our calculations. While our method does not find a simple ``closed form" for the density of the invariant measure (if one even exists), our work leads us to some new conjectures about the behavior of the density at certain points.

In addition to this, we  detail all admissible strings of digits in the Hurwitz expansion. This may be of independent interest.
\end{abstract}

\section{Introduction}

Given a transformation of a subset of $\mathbb{R}^n$ to itself, one common question is whether or not there exists an absolutely continuous (with respect to Lebesgue) invariant measure (or a.c.i.m), and if it exists, to find a simple closed form of it. This is a very common concern in the study of fibred systems, the symbolic expansions of numbers. The transformation $Tx=bx\pmod{1}$ on $[0,1)$ associated to the base-$b$ transformation has Lebesgue itself as its invariant measure. The transformation $Tx=\beta x\pmod{1}$ on $[0,1)$ associated to $\beta$-expansions ($\beta>1$) has the Parry measure as its invariant measure. The transformation $Tx=1/x\pmod{1}$ on $[0,1)$ associated to the regular continued fraction expansion has the Gauss measure $\mu(A)=\int_A \frac{dx}{\log 2(1+x)}$ as its invariant measure. (See \cite{DK}.)

Of particular interest to us are the invariant measures of the many continued fraction (CF) variants. Explicit closed forms are known for the invariant measure of the backwards CF, even CF, odd CF, Rosen CF, and some, but not all, of Nakada's $\alpha$-CFs \cite{BKS,Masarotto}. These are all one-dimensional real CF expansions. The simplest complex CF expansion is the Hurwitz complex CF (see Section \ref{sec:fundamentals} for definitions).\footnote{Brothers Adolf and Julius Hurwitz each have their own complex continued fraction expansion \cite{OS}. We will be considering the expansion investigated by Adolf Hurwitz.} While the Hurwitz complex CF is well-studied \cite{BG,DN,Hensley,NakadaKuzmin,Schweiger00}, little is yet known about the corresponding invariant measure. It is known (see Theorem \ref{thm:Hensley}) that the density $h$ of the invariant measure is piece-wise real-analytic with $12$ pieces of analyticity and it is known that it satisfies certain symmetries, but that is all. Some of Doug Hensley's computational work on the invariant density $h$ can be found at his website: \url{http://www.math.tamu.edu/~dhensley/}.

In Section \ref{sec:method} of this paper, we will describe a method for numerically approximating a truncated Taylor expansion of the invariant density $h$ around a given point. In Section \ref{sec:results}, we describe the results of our computation. The precision of our calculations relied on a variable $k$, which determines the quality of the numerical approximation. We found that increasing $k$ by one roughly reduced the error in our computation of the coefficients by a multiplicative factor of around $.57$ (see Table \ref{tab:V11abc}), which suggests that as we increase $k$ linearly, our numerical approximations converge exponentially. At the limit of computation on our personal computers, this suggests we were able to compute the coefficients to an accuracy of around $\pm 10^{-3}$. As an example of our calculations, we approximated $h$ near the point $(-.5,-.5)$ (using the usual identification of $\mathbb{C}$ with $\mathbb{R}^2$) with the following function (again, see Table \ref{tab:V11abc}):
\begin{align*}
&0.7149+0.3411(x+.5)^2+0.3411 (y+.5)^2+0.0875(x+.5)^4+0.0875(y+.5)^4\\
&\qquad +0.4974 (x+.5)^2(y+.5)^2+0.0180(x+.5)^6+0.0180(y+.5)^6 \\
&\qquad+0.2652(x+.5)^4(y+.5)^2+0.2652(x+.5)^2(y+.5)^4
\end{align*}
Mathematica files including the implementation of our calculations can be found at \url{http://github.com/JVandehey/HurwitzInvariantMeasure}.

\begin{figure}[h!]
\caption{A plot of the density function $h(z)$ of the invariant measure}\label{fig:1}
\includegraphics[width=\textwidth]{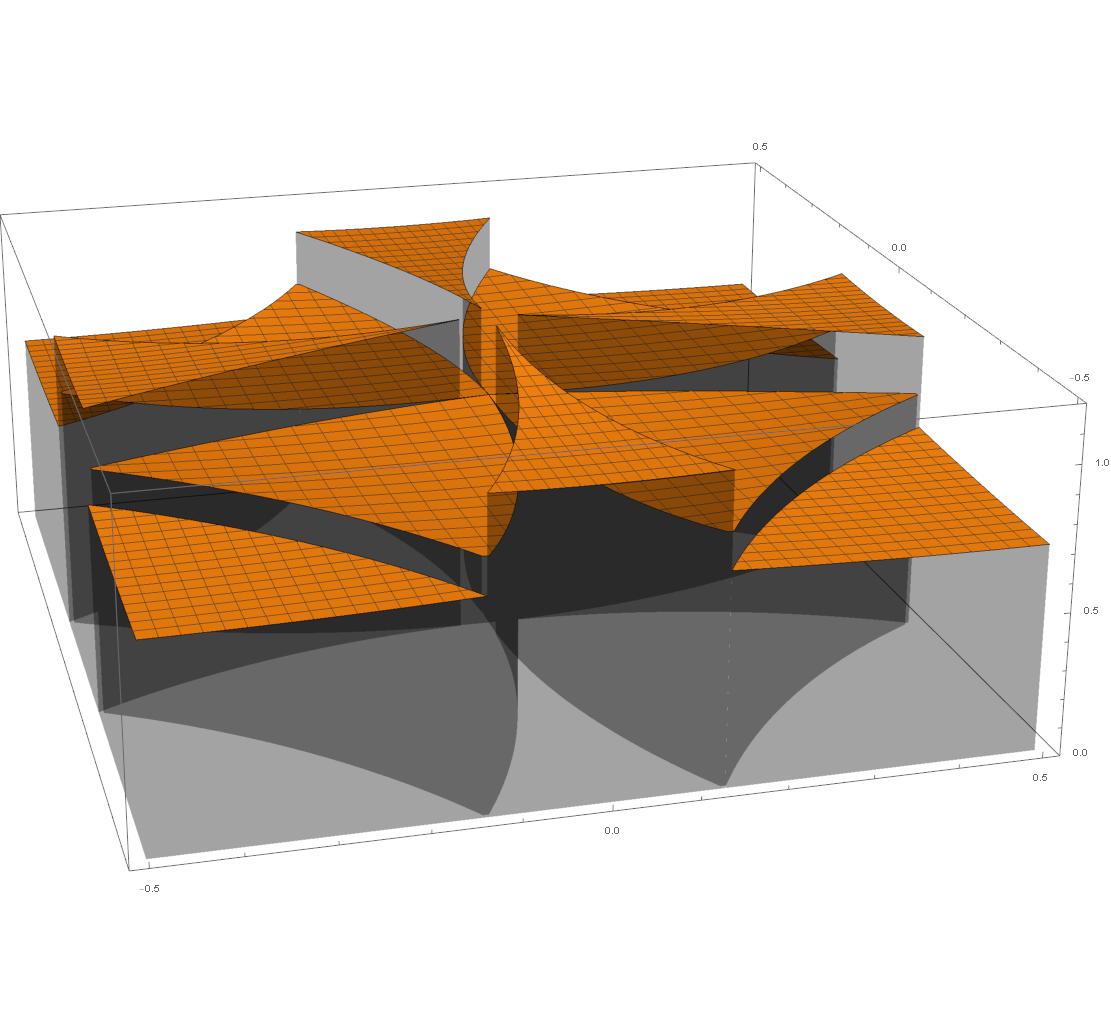}
\end{figure}

\begin{figure}[h!]
\caption{A zoomed-in plot of three regions of the density function $h(z)$ of the invariant measure}\label{fig:2}
\includegraphics[width=\textwidth]{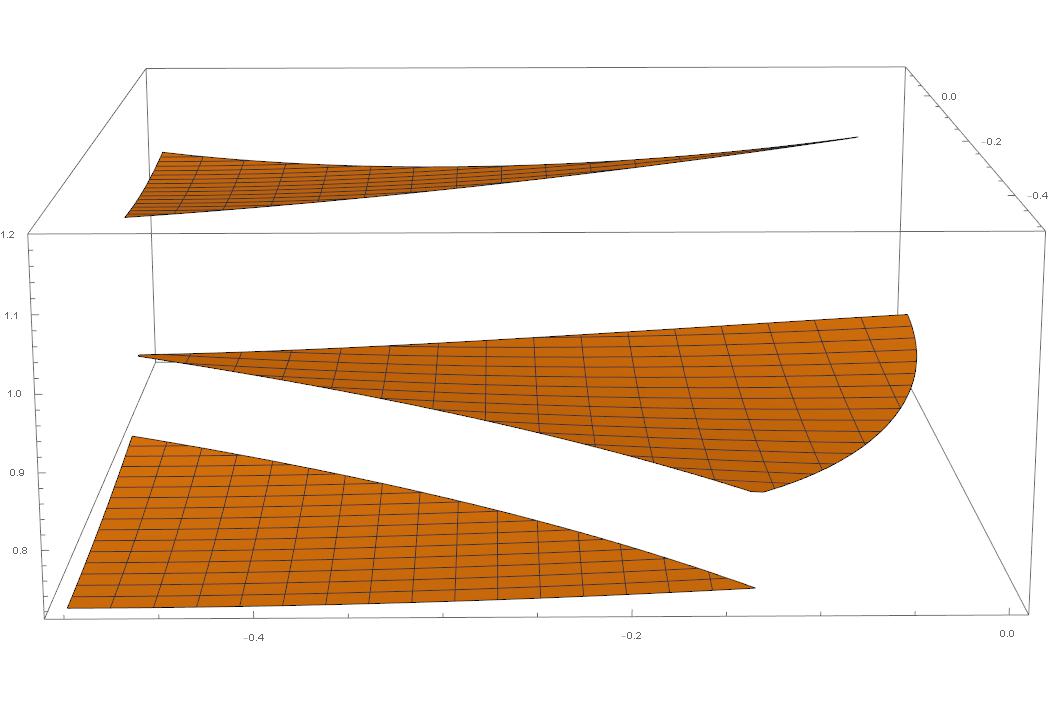}
\end{figure}

Figure \ref{fig:1} shows the complete invariant density produced using our data, with all 12 regions of analyticity and $90$ degree rotational symmetry around the origin clearly visible. It is  difficult to see in this picture, but all 12 regions show some curvature to them. We also produced Figure \ref{fig:2}, which shows an enlarged image of 3 of the regions of analyticity, with more of the curvature now visible.

While this information does not give a simple, closed-form expression for the invariant density, our calculations did suggest the following conjecture:

\begin{conj}
The Taylor series of the invariant density $h$ is an even function in both the $x$ and $y$ coordinates around the points $\pm .5\pm .5\ii$, $\pm .5$, and $\pm .5\ii$. In other words, if $x_0+y_0\ii$ is any of the above points, then
\[
\left.\frac{\partial^{m+n}}{\partial^m x\partial^n y} h(x+y\ii) \right|_{x+y\ii=x_0+y_0\ii}
\]
is $0$ whenever $m$ or $n$ is odd.
\end{conj}

We are being somewhat informal in talking about the Taylor series. In particular, when we discuss the Taylor series of $h$ at a point, we mean the Taylor series of the corresponding analytic piece that the point lies in. Also, each analytic piece of $h$ could be extended to an analytic function on a larger set (this is a consequence of a proof in \cite{Hensley}), and thus it is reasonable to discuss the Taylor series at points on the boundary of the domain of $h$.

In Section \ref{sec:othermethods}, we briefly describe some other methods we used to try and compute the invariant density and why they didn't work. One method lead us to the question of what admissible sequences of digits can appear in the Hurwitz complex CF expansion. We answer that in Section \ref{sec:admissible}. One consequence of our work is the following result.

\begin{thm}\label{thm:admissible}
If one allows the digits of a Hurwitz CF expansion to include Gaussian integers and marked Gaussian integers, then admissible strings of digits can be completely determined by the length-$2$ admissible strings of digits. That is, if $[a_1,a_2,\dots,a_n]$ is an admissible string, then the possible values for $a_{n+1}$ such that $[a_1,a_2,\dots,a_n,a_{n+1}]$ is also admissible are determined solely by the value of $a_n$.
\end{thm}

We note that marked numbers, such as $(1+\ii)'$, are distinct from regular numbers, such as $1+\ii$, only symbolically. They take the same numerical value.

\section{Fundamentals of Hurwitz Continued Fractions}\label{sec:fundamentals}

Let $K\subset\mathbb{C}$ denote the set \[ K:= \{x+y \ii:x,y\in[-1/2,1/2)\}. \]
It is clear that the shifts of this set by Gaussian integers tesselate the complex plane. By this we mean that $\bigcup_{\alpha\in\mathbb{Z}[\ii]} (K+\alpha)=\mathbb{C}$ and for any two distinct $\alpha,\alpha'\in \mathbb{Z}[\ii]$, we have \[ (K+\alpha)\cap (K+\alpha')=\emptyset.\]

As such, to each element $z\in \mathbb{C}$ we may define a unique element $[z]\in \mathbb{Z}[\ii]$ such that $z\in K+[z]$. This function $[\cdot]$ we may interpret as the nearest integer to $z$ and can write explicitly in terms of the real-valued nearest-integer function by $[x+y \ii]=[x]+\ii[y]$.

From this we may define a Gauss map $T:K\to K$ given by 
\[
T(z) = \begin{cases} \dfrac{1}{z}-\left[ \dfrac{1}{z} \right], & \text{ if }z\neq 0,\\
0, & \text{ if }z=0.\end{cases}
\]
For a given point $z\in K$ we define (possibly finite, possibly infinite) sequences $z_n$, $a_n$, defined inductively by $z_0=z$ and then, provided $z_n\neq 0$, we let $a_{n+1}=[1/z_n]$ and $z_{n+1}=(1/z_n)-a_{n+1}=T(z_n)$. If $z_n=0$, then we terminate the sequences at $n$. From this we see that
\[
z=\dfrac{1}{a_1+z_1}=\dfrac{1}{a_1+\dfrac{1}{a_2+z_2}}=\dfrac{1}{a_1+\dfrac{1}{a_2+\dfrac{1}{a_3+z_3}}}.
\]

We would like to prove the meaningfulness and convergence of the infinite continued fraction expansion,
\[
z=\dfrac{1}{a_1+\dfrac{1}{a_2+\dfrac{1}{a_3+\dots}}}.
\]
The convergence is immediate if the sequence is finite (which happens if and only if $z\in\mathbb{Q}[\ii]$). To indicate why convergence happens in the infinite case, we must consider the convergents. We let $p_{-1}=q_0=1$ and $p_0=q_{-1}=0$, and then let for $n\ge 1$, let $p_n=a_np_{n-1}+p_{n-2}$ and $q_n=a_nq_{n-1}-q_{n-2}$. It is an easy induction proof to show that 
\[
\frac{p_n}{q_n} =\dfrac{1}{a_1+\dfrac{1}{a_2+\dfrac{1}{a_3+\dots+\frac{1}{a_n}}}},
\]
and that this fraction is in lowest terms. Moreover, with slightly more difficulty, one can show that
\[
\left| z-\frac{p_n}{q_n} \right| < \frac{2\sqrt{2}|z_n|}{|q_n|^2}
\]
and that $|q_{n+2}/q_n|\ge 3/2$, proving that the convergents $p_n/q_n$ do indeed converge to $z$ as desired. (For details, see \cite[Ch.~5]{Hensley})

Since one can write $a_{n+1}(z)=[1/T^n z]$, we also have that $a_n(Tz)=[1/T^n z]$, and hence $T$ acts as a forward-shift on the string of digits, i.e., if we connect a point $z\in K$ to its string of digits $[a_1,a_2,a_3,\dots]$, then  \begin{equation}\label{eq:Tshift} T([a_1,a_2,a_3,\dots])=[a_2,a_3,a_4,\dots].\end{equation}

Particularly relevant for our purposes, we know that there exists a probability measure $\mu$ on $K$, absolutely continuous with respect to Lebesgue, that is invariant---that is, $\mu(T^{-1}A)=\mu(A)$---and ergodic---that is, if $T^{-1}A=A$, then $\mu(A)=0$ or $\mu(A)=1$---with respect to the transformation $T$. This was first shown by Nakada \cite{NakadaKuzmin}. As a consequence of further investigations into the map $T$, Schweiger \cite{Schweiger00,Schweiger11} showed that $\mu$ is piecewise Lipschitz with respect to a certain partition of $K$. Hensley \cite{Hensley} further proves the following result.

\begin{thm}\label{thm:Hensley}
The measure $\mu$ has a density function $\rho$---that is, $\mu(A) = \int_A \rho \ d\lambda(z)$ with $\lambda$ being Lebesgue measure on $\mathbb{C}$---where $\rho$ is continuous except perhaps along the arcs $|z\pm 1|=1$, $|z\pm \ii|=1$, and $|z\pm 1\pm \ii|=1$. It is moreover real-analytic on each of the 12 open regions that these arcs partition $K$ into. $\mu$ also obeys the following symmetries: $\mu(	\ii A)=\mu(A)$ and $\mu(\overline{A})=\mu(A)$.
\end{thm}

In order to better understand the measure $\mu$, we introduce the natural extension. Let $\mathbb{C}^*=\mathbb{C}\cup \{\infty\}$

We extend $T$ to $\hat{T}:K\times\mathbb{C}^*\to K\times \mathbb{C}^*$, that acts by
\[
\hat{T}(z,w)=\begin{cases}\left( \dfrac{1}{z}-\left[\dfrac{1}{z}\right], \dfrac{1}{\left[\frac{1}{z}\right]+w}\right), & z\neq 0,\\
\left( z,w\right), & z=0.\end{cases}
\]
We then define $\hat{K}\subset (\mathbb{C}^*)^2$ by
\[
\hat{K}=\overline{\left\{\hat{T}^i (z,0): z\in K, i\in \mathbb{Z}_{\ge 0}\right\}}.
\]
Whereas typically we use the overline to refer to taking conjugates, in this case we use it to define the closure. Note that it can be shown further (see again \cite[Ch.~5]{Hensley}) that $\hat{K}\subset B^2$ where $B$ is the unit ball around the origin in $\mathbb{C}$.

Similar to \eqref{eq:Tshift}, we have that
\[
\hat{T}^i \left( \dfrac{1}{a_1+\dfrac{1}{a_2+\dfrac{1}{a_3+\dots}}},0\right) = \left( \dfrac{1}{a_{i+1}+\dfrac{1}{a_{i+2}+\dfrac{1}{a_{i+3}+\dots}}}, \dfrac{1}{a_{i}+\dfrac{1}{a_{i-1}+\dots + \dfrac{1}{a_1}}}\right).
\]
It is true, although not immediate, that by ignoring an appropriate measure-zero set, $\hat{T}$ acts as a bijection from $\hat{K}$ to itself, as was shown in Ei, Ito, Nakada, and Natsui \cite{EINN}.

Let $V\subset \mathbb{C}$ denote the projection of $\hat{K}$ onto its second coordinate. It is not the case that $\hat{K}=K\times V$ (up to some zero-measure set); however, it is thankfully not much more complicated than that. 

\begin{figure}[h!]
\caption{Some pictures of $K_{k,\ell}$ and $V_{k,\ell}$. Note that the picture of $V_{3,1}$ is scaled differently than $V_{1,1}$ and $V_{2,1}$.}
\centering
\begin{subfigure}[b]{0.3\textwidth}
\includegraphics[width=\textwidth]{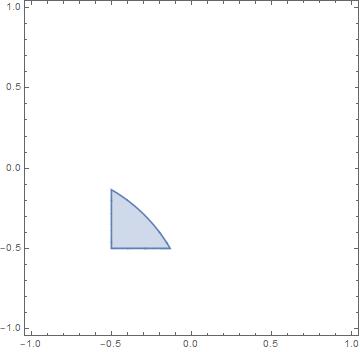}
\caption{The region $K_{1,1}$}
\label{fig:K}
\end{subfigure}
\qquad
\begin{subfigure}[b]{0.47\textwidth}
\includegraphics[width=\textwidth]{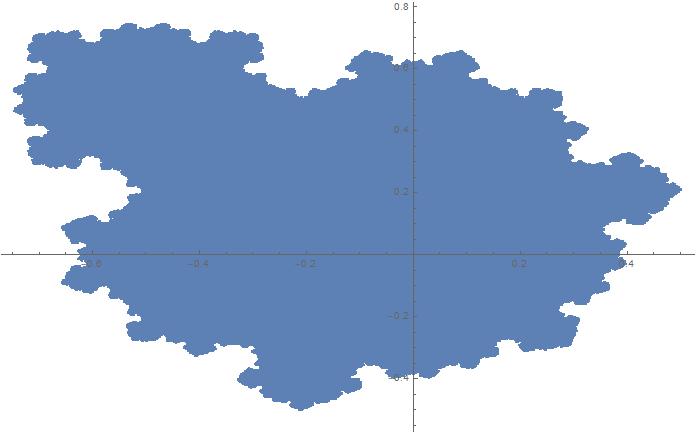}
\caption{The region $V_{1,1}$}
\label{fig:Kinv}
\end{subfigure}

\

\begin{subfigure}[b]{0.3\textwidth}
\includegraphics[width=\textwidth]{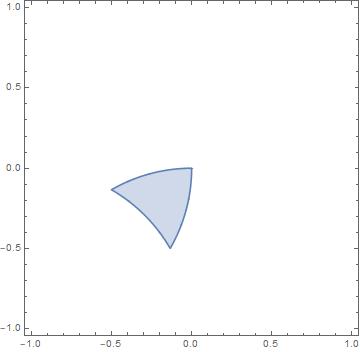}
\caption{The region $K_{2,1}$}
\label{fig:K}
\end{subfigure}
\qquad
\begin{subfigure}[b]{0.47\textwidth}
\includegraphics[width=\textwidth]{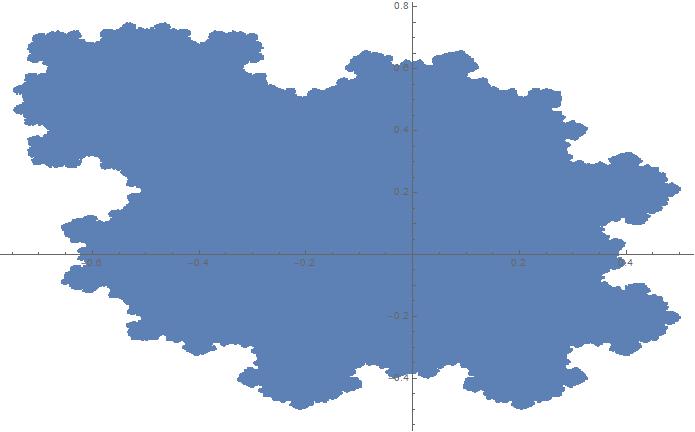}
\caption{The region $V_{2,1}$}
\label{fig:Kinv}
\end{subfigure}

\

\begin{subfigure}[b]{0.3\textwidth}
\includegraphics[width=\textwidth]{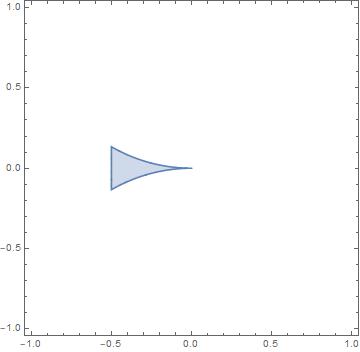}
\caption{The region $K_{3,1}$}
\label{fig:K}
\end{subfigure}
\qquad
\begin{subfigure}[b]{0.47\textwidth}
\includegraphics[width=\textwidth]{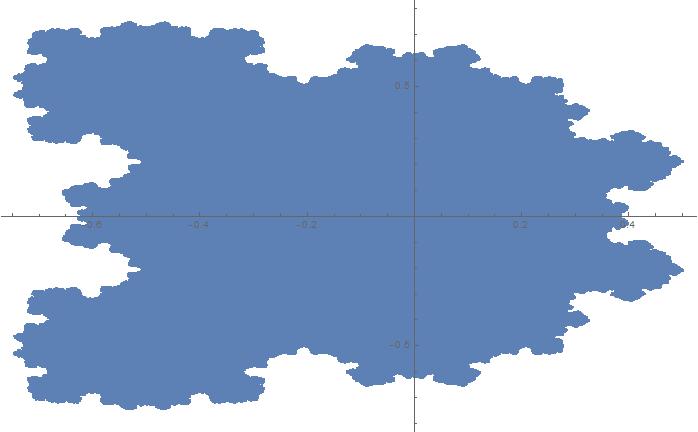}
\caption{The region $V_{3,1}$}
\label{fig:Kinv}
\end{subfigure}

\caption{}
\label{fig:KVdecomp}

\end{figure}

Let us define the 12 open regions of $K$ more precisely.\footnote{Although we use somewhat similar notation to the Ei, et al, paper, our definitions are distinct and should not be mistaken for one another.} Let $K_{k,\ell}$ for $k\in\{1,2,3\}$, $\ell\in \{1,2,3,4\}$, be given by
\begin{align}
K_{1,\ell}&= \ii^{\ell-1} \cdot \{z\in K: |z+1+\ii|<1\}\notag\\
K_{2,\ell}&= \ii^{\ell-1} \cdot \{z\in K: |z+1+\ii|>1,|z+1|<1,|z+\ii|<1\}\label{eq:Kkl}\\
K_{3,\ell}&= \ii^{\ell-1} \cdot \{z\in K: |z+1|<1,|z\pm \ii |>1\}.\notag
\end{align}
Note that these cover $K$ up to a zero-measure set. (See the left-hand side of Figure \ref{fig:KVdecomp}.) 

Now to each $K_{k,\ell}$ there exists a set $V_{k,\ell}\subset V$ such that $\hat{K}= \bigcup_{k,\ell} K_{k,\ell}\times V_{k,\ell}$, again up to a zero-measure set. (For a depiction of some of the $V_{k,\ell}$, see the right-hand side of Figure \ref{fig:KVdecomp}.) We note several symmetries are apparent from the pictures and are in fact as we see them. For example, $V_{1,1}$ is symmetric along the line $b=-a$, if we think of complex points as being written as $a+b\ii$. Likewise $V_{2,1}$ is also symmetric along $b=-a$ and $V_{3,1}$ is symmetric along $b=0$.

 Ei, et al, \cite{EINN} study the sets $V_{k,\ell}$ very closely and in particular prove a number of useful facts. In particular, each $V_{k,\ell}$ has positive Lebesgue measure with path-connected and simply connected interior and a boundary given by a Jordan curve, which appears to be fractal in nature.

\subsection{Properties of the invariant measure}

These properties allow us to nicely define an absolutely continuous invariant measure $\hat{\mu}$ for $\hat{T}$ on $\hat{K}$. In particular, we will show that if we define $\hat{\mu}$ by 
\[
\hat{\mu}(E) = \int_E \frac{1}{|zw+1|^4} \ d\lambda(z) \ d\lambda(w), \ \text{where } (z,w)\in E
\]
for any measurable set $E\subset \hat{K}$, then this is an invariant measure. That this integral is actually integrable follows from the properties in the previous paragraph. That this measure is invariant is fairly straight-forward. Since $\hat{T}$ is a bijection, it suffices to show $\hat{\mu}(\hat{T}E)=\hat{\mu}(E)$, for any measurable $E\subset \hat{K}$. Moreover, since $\hat{\mu}$ is countably additive and there are countably many digits, it suffices to consider a set $E$ such that for each $(z,w)\in E$, the first continued fraction digit of $z$ is the same. Thus $\hat{T}$ acts by
\[
\hat{T}(z,w)= \left( \frac{1}{z} - a, \frac{1}{a+w}\right),
\]
where $a=a_1(z)$. Thus,
\begin{align*}
 \int_{\hat{T}E} \frac{1}{|zw+1|^4}  \ d\lambda(z) \ d\lambda(w)  &=\int_E \frac{1}{\left| \left( \frac{1}{z}-a\right) \left( \frac{1}{a+w}\right)+1 \right|^4} \frac{d\lambda(z)}{|z|^4} \ \frac{d\lambda(w)}{|a+w|^4}\\
&= \int_E \frac{1}{\left| (1-za)+z(a+w)\right|^4}  \ d\lambda(z) \ d\lambda(w)\\
&= \int_E \frac{1}{|zw+1|^4}  \ d\lambda(z) \ d\lambda(w).
\end{align*}
We note briefly that when we do the change of variables, we get $d\lambda(z)/|z|^4$ instead of $d\lambda(z)/|z|^2$ because  $d\lambda(z)$ is the differential with respect to Lebesgue measure and since this is seen as a transformation on the complex plane, areas are shrunk by the square of the derivative. 

So this does define an invariant measure for $\hat{T}$ on $\hat{K}$.

This gives us a way to define $\mu$ by projecting from $\hat{\mu}$. Let $\pi_K:\hat{K}\to K$ be the projection onto the first coordinate. (Since $\hat{K}$ was defined using a closure and $K$ does not contain its entire boundary, this projection will not be defined for a Lebesgue-measure-zero subset of $\hat{K}$, but this will ultimately not matter.) 

Briefly, we will show that $\pi_K^{-1}T^{-1}=\hat{T}^{-1}\pi_K^{-1}$. Namely, $(z,w)\in \pi_K^{-1} T^{-1} A$ if and only if $z\in T^{-1} A$ if and only if $Tz\in A$. Moreover, $(z,w)\in \hat{T}^{-1}\pi_K^{-1} A$ if and only if $\hat{T}(z,w)\in \pi_K^{-1} A$ if and only if $Tz \in \pi_K^{-1} A$. 

And thus if we define a measure $\mu=\hat{\mu}\circ \pi_K^{-1}$ then for any measurable set $A\subset K$, we have
\[
\mu(T^{-1} A)= \hat{\mu}(\pi_K^{-1}T^{-1}A )=\hat{\mu}(\hat{T}^{-1}\pi_K^{-1} A) = \hat{\mu}(\pi_K^{-1} A)= \mu(A),
\]
so this measure is invariant and one can also see that it is absolutely continuous by the absolute continuity of $\hat{\mu}$. So this measure must be the same $\mu$ as before. (We are eliding over a small step here: to show this measure and the previous one were the same, we must invoke that we already knew the previous measure was ergodic and use the ergodic decomposition theorem.)

Thus, if we only knew the shapes $V_{k,\ell}$ perfectly (and perfectly how to integrate over them), then the problem of understanding $\mu$ would be trivialized. While we cannot yet do that, we can use the simple form of $\hat{\mu}$ to help us estimate $\mu$ better. 

In particular, suppose $h$ is the density of $\mu$ (with respect to Lebesgue measure). Then we can express $h$ by
\[
h(z)= \int_{V_{k,\ell}} \frac{1}{|zw+1|^4} \ dw, \qquad \text{ if }z\in K_{k,\ell}
\]
We may express this as a function $h(x,y)$ where $z=x+iy$ and since it is real-analytic, we may express its derivatives via
\begin{equation}\label{eq:hderivative}
h_{m,n}(x,y):=\frac{\partial^{m+n}}{\partial x^m \partial y^n} h(x,y) = \int_{a+b\ii\in V_{k,\ell}}\frac{\partial^{m+n}}{\partial x^m \partial y^n} \frac{1}{((ax-by+1)^2+(ay+bx)^2)^2} \ da \ db.
\end{equation}
This equation comes from taking $w=a+b\ii$.

When the choice of $k,\ell$ is not clear, we will use $h_{m,n}^{k,\ell}(x,y)$. In general, we will assume that $x+y\ii\in K_{k,\ell}$, but in at least one case, we will not make this assumption.

Equation \eqref{eq:hderivative} will be the primary tool we use to calculate the Taylor series for $h$ at appropriate points.

Finally, we conclude this section with a discussion of the Perron-Frobenius operator. This operator $P$ is defined in this case by
\[
\int_A Pf(z) d\lambda(z) = \int_{T^{-1}A} f(z) d\lambda(z), \qquad \text{for all measurable }A\subset K.
\]
If $f=h$, where $h$ is the density of the measure $\mu$, then clearly $Pf=f$. On the other hand, if $Pf=f$, then by the above equation, the measure $\mu_f$ given by $\mu_f(A)=\int_A f d\lambda$ is $T$-invariant. By the ergodic decomposition theorem, we know that $\mu_f$ can be, loosely speaking, written as a linear combination of ergodic measures. We also know that any two ergodic measures are mutually singular---that is, the support of one ergodic measure is of measure zero with respect to any other ergodic measure---so since we already know that $\mu$ is ergodic and has $K$ as its support, the support of any other ergodic measure must have zero Lebesgue measure. Since $\mu_f$ is continuous with respect to Lebesgue, it must therefore be a constant times $\mu$, and so, up to a constant, $f$ is equal to $h$ almost everywhere.

Moreover, it can be shown that
\[
Pf(z) = \sum_{z'\in T^{-1}z} \frac{f(z')}{|T'(z')|},
\]
where here the denominator represents the Jacobian of $T$ seen here as two-dimensional real-valued map. This gives us a functional equation for the invariant density $h$, which could also be used to calculate it. In particular,
\[
h(z)= \sum_{(a+z)^{-1}\in T^{-1}z} \frac{1}{|a+z|^4}h\left(\dfrac{1}{a+z}\right), \ z\in K.
\]

\section{Admissible strings of digits}\label{sec:admissible}

In Figure \ref{fig:K}, we see the region $K$. When we try to calculate $z_1$, the first thing we must do is invert, i.e., apply the map $z \mapsto 1/z$. This map applied to the region $K$ gives 
Figure \ref{fig:Kinv}. We have divided the map up into tiles centered around each Gaussian integer, each of which corresponds to a unique possibility of $a_1$, with the correspondence being the obvious one: each tile corresponds to the Gaussian integer it is centered at. If there is any intersection between the tile and the inversion of $K$, we know that that digit is a possibility for $a_1$. In particular, $a_1$ can be any element of the set 
\[
G=\mathbb{Z}[\ii]\setminus\{0,1,-1,\ii,-\ii\}.
\]

\begin{figure}[h!]
\centering
\begin{subfigure}[b]{0.25\textwidth}
\includegraphics[width=\textwidth]{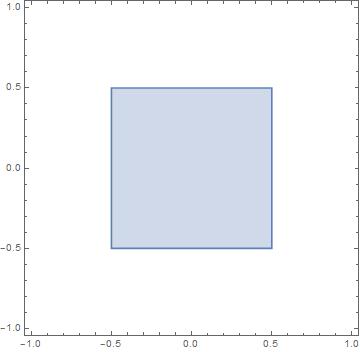}
\caption{The region $K$}
\label{fig:K}
\end{subfigure}
\qquad
\begin{subfigure}[b]{0.25\textwidth}
\includegraphics[width=\textwidth]{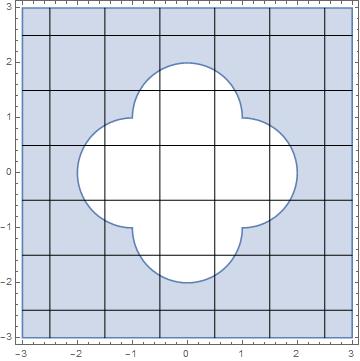}
\caption{The region $K$ inverted}
\label{fig:Kinv}
\end{subfigure}

\caption{The initial picture}
\label{fig:Kinitial}

\end{figure}

This picture is extremely important. Although $z=z_0$ could be anything in $K$, this is not true of $z_1$. The possible values of $z_1$ are dependent on $a_1$. For example, if $a_1=2$ (which is possible because in Figure \ref{fig:Kinv} we see that the shaded region intersects the tile centered at $2$) then we cannot have that $z_1=-1/2$, because when we shift this tile back to be centered at the origin the point $-1/2$ is not included in the shaded region. 

Extending this analysis slightly, if $a_1=2$, then $z_1$ cannot satisfy $|z_1+1|<1$, and then, by inverting the corresponding circle, we see that $\Re(1/z_1)\ge -1/2$, and thus we cannot have any digit $a_2$ for which the real part is negative. This corresponds to part of the first line in Table \ref{table:allowabledigits}.

\begin{figure}[h!]
\centering
\begin{subfigure}[b]{0.25\textwidth}
\includegraphics[width=\textwidth]{Digits1}
\end{subfigure}
\qquad
\begin{subfigure}[b]{0.25\textwidth}
\includegraphics[width=\textwidth]{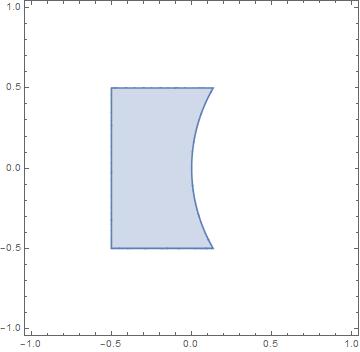}
\end{subfigure}

\

\begin{subfigure}[b]{0.25\textwidth}
\includegraphics[width=\textwidth]{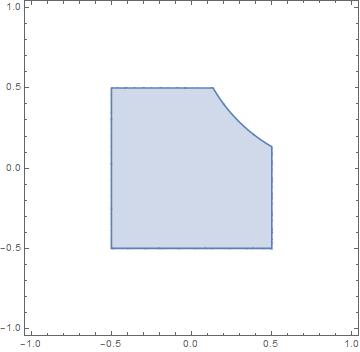}
\end{subfigure}
\qquad
\begin{subfigure}[b]{0.25\textwidth}
\includegraphics[width=\textwidth]{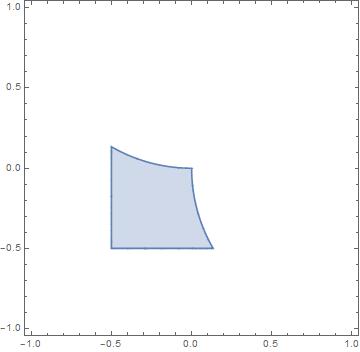}
\end{subfigure}

\caption{The four possible new regions (up to rotation)}

\label{fig:Krotation}

\end{figure}

However, even though the possible regions $z_1$ can exist in depend on $a_1$, they are not extremely complex. There are 13 possible regions: the full region, and 12 subregions. The 12 subregions can be broken down further into 3 examples, with the remaining 9 obtained by rotating these by $\pi/2$, $\pi$, or $3\pi/2$. The full region, as well as the 3 initial subregions, are given in Figure \ref{fig:Krotation}. We note that these regions are obtained by taking unions of the various $K_{k,\ell}$ regions defined in \eqref{eq:Kkl}, up to a finite measure set along their boundaries.

\begin{figure}
\centering
\begin{subfigure}[b]{0.25\textwidth}
\includegraphics[width=\textwidth]{Digits2}
\caption{A subregion...}
\label{fig:D2}
\end{subfigure}
\qquad
\begin{subfigure}[b]{0.25\textwidth}
\includegraphics[width=\textwidth]{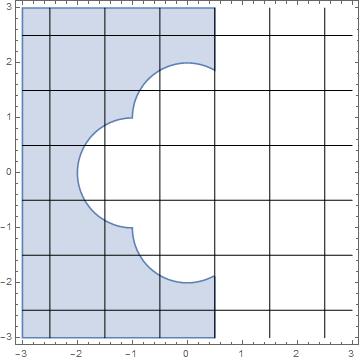}
\caption{...and its inverse}
\label{fig:D2inv}
\end{subfigure}

\

\begin{subfigure}[b]{0.25\textwidth}
\includegraphics[width=\textwidth]{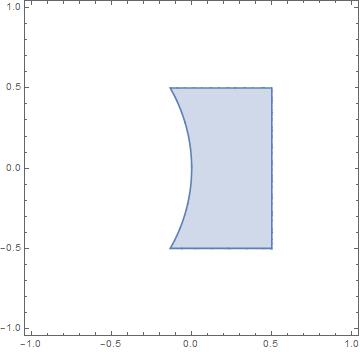}
\caption{A subregion...}
\label{fig:D3}
\end{subfigure}
\qquad
\begin{subfigure}[b]{0.25\textwidth}
\includegraphics[width=\textwidth]{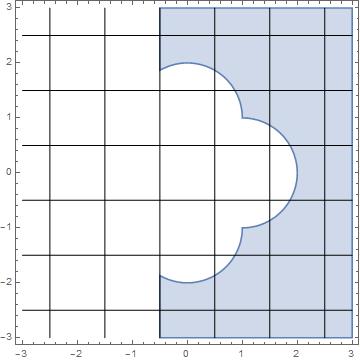}
\caption{...and its inverse}
\label{fig:D3inv}
\end{subfigure}

\

\begin{subfigure}[b]{0.25\textwidth}
\includegraphics[width=\textwidth]{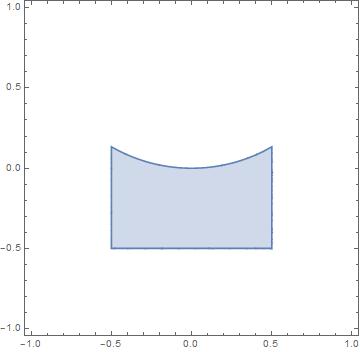}
\caption{A subregion...}
\label{fig:D4}
\end{subfigure}
\qquad
\begin{subfigure}[b]{0.25\textwidth}
\includegraphics[width=\textwidth]{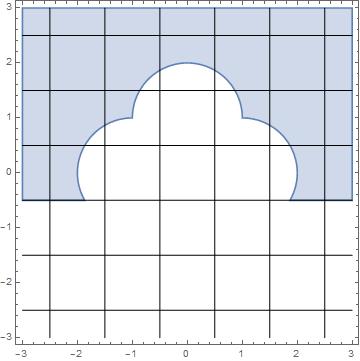}
\caption{...and its inverse}
\label{fig:D4inv}
\end{subfigure}

\

\begin{subfigure}[b]{0.25\textwidth}
\includegraphics[width=\textwidth]{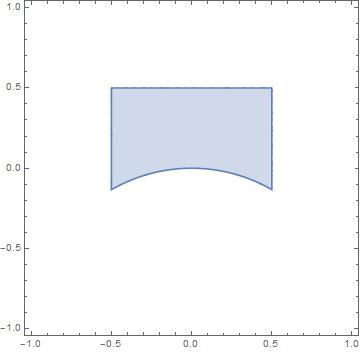}
\caption{A subregion...}
\label{fig:D5}
\end{subfigure}
\qquad
\begin{subfigure}[b]{0.25\textwidth}
\includegraphics[width=\textwidth]{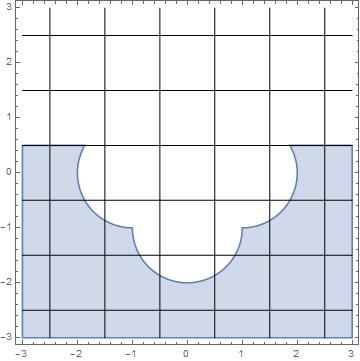}
\caption{...and its inverse}
\label{fig:D5inv}
\end{subfigure}

\

\caption{The first subregion, its rotations, and their inversions}

\label{fig:firstsubregion}

\end{figure}

\begin{figure}
\centering
\begin{subfigure}[b]{0.25\textwidth}
\includegraphics[width=\textwidth]{Digits6}
\caption{A subregion...}
\label{fig:D6}
\end{subfigure}
\qquad
\begin{subfigure}[b]{0.25\textwidth}
\includegraphics[width=\textwidth]{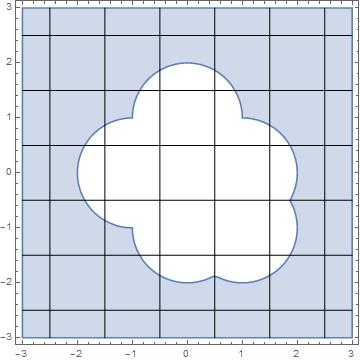}
\caption{...and its inverse}
\label{fig:D6inv}
\end{subfigure}

\

\begin{subfigure}[b]{0.25\textwidth}
\includegraphics[width=\textwidth]{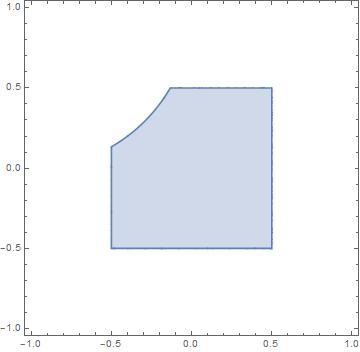}
\caption{A subregion...}
\label{fig:D7}
\end{subfigure}
\qquad
\begin{subfigure}[b]{0.25\textwidth}
\includegraphics[width=\textwidth]{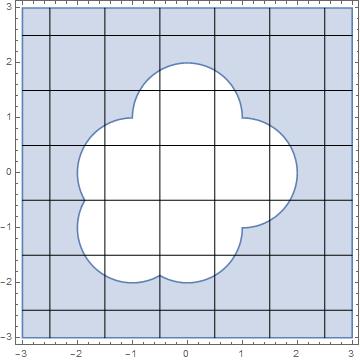}
\caption{...and its inverse}
\label{fig:D7inv}
\end{subfigure}

\

\begin{subfigure}[b]{0.25\textwidth}
\includegraphics[width=\textwidth]{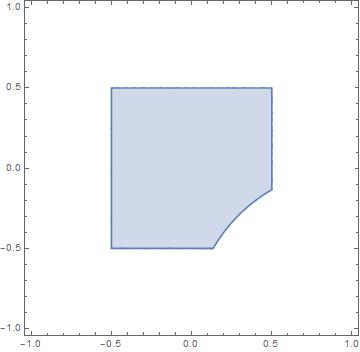}
\caption{A subregion...}
\label{fig:D8}
\end{subfigure}
\qquad
\begin{subfigure}[b]{0.25\textwidth}
\includegraphics[width=\textwidth]{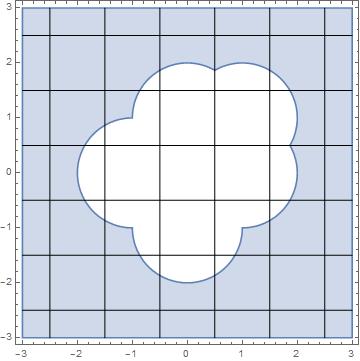}
\caption{...and its inverse}
\label{fig:D8inv}
\end{subfigure}

\

\begin{subfigure}[b]{0.25\textwidth}
\includegraphics[width=\textwidth]{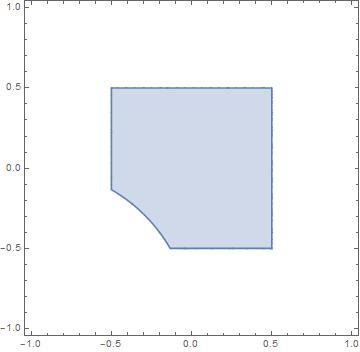}
\caption{A subregion...}
\label{fig:D9}
\end{subfigure}
\qquad
\begin{subfigure}[b]{0.25\textwidth}
\includegraphics[width=\textwidth]{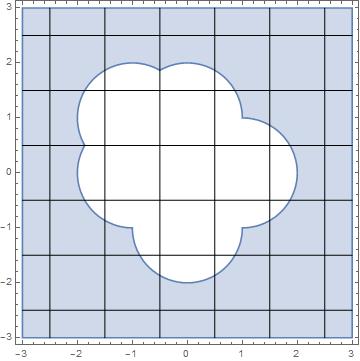}
\caption{...and its inverse}
\label{fig:D9inv}
\end{subfigure}

\

\caption{The second subregion, its rotations, and their inversions}

\label{fig:secondsubregion}

\end{figure}

\begin{figure}
\centering
\begin{subfigure}[b]{0.25\textwidth}
\includegraphics[width=\textwidth]{Digits10}
\caption{A subregion...}
\label{fig:D10}
\end{subfigure}
\qquad
\begin{subfigure}[b]{0.25\textwidth}
\includegraphics[width=\textwidth]{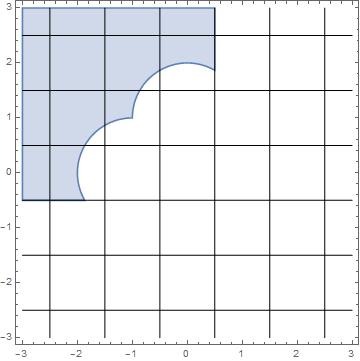}
\caption{...and its inverse}
\label{fig:D10inv}
\end{subfigure}

\

\begin{subfigure}[b]{0.25\textwidth}
\includegraphics[width=\textwidth]{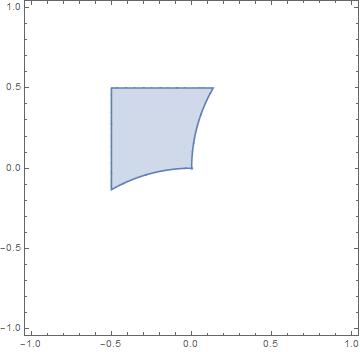}
\caption{A subregion...}
\label{fig:D11}
\end{subfigure}
\qquad
\begin{subfigure}[b]{0.25\textwidth}
\includegraphics[width=\textwidth]{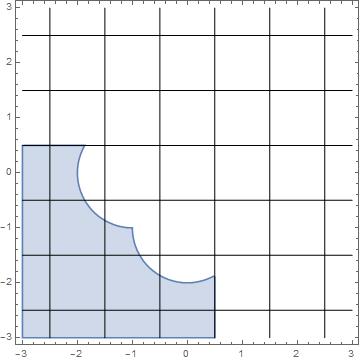}
\caption{...and its inverse}
\label{fig:D11inv}
\end{subfigure}

\

\begin{subfigure}[b]{0.25\textwidth}
\includegraphics[width=\textwidth]{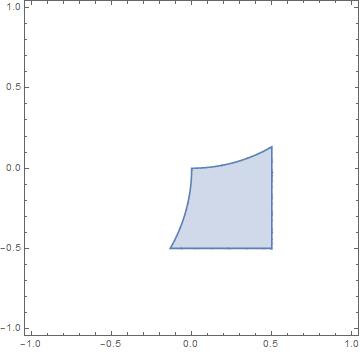}
\caption{A subregion...}
\label{fig:D12}
\end{subfigure}
\qquad
\begin{subfigure}[b]{0.25\textwidth}
\includegraphics[width=\textwidth]{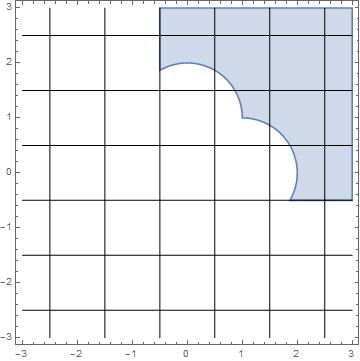}
\caption{...and its inverse}
\label{fig:D12inv}
\end{subfigure}

\

\begin{subfigure}[b]{0.25\textwidth}
\includegraphics[width=\textwidth]{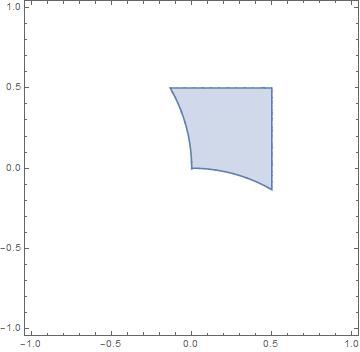}
\caption{A subregion...}
\label{fig:D13}
\end{subfigure}
\qquad
\begin{subfigure}[b]{0.25\textwidth}
\includegraphics[width=\textwidth]{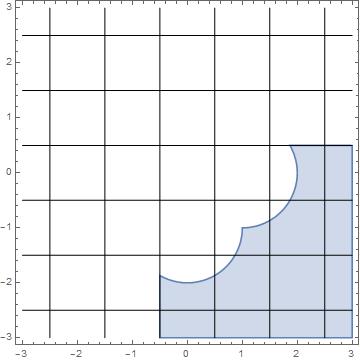}
\caption{...and its inverse}
\label{fig:D13inv}
\end{subfigure}

\

\caption{The third subregion, its rotations, and their inversions}

\label{fig:thirdsubregion}

\end{figure}

Now, in order to figure out the possibilities for $a_2$ (depending on $a_1$, of course), we must take each of these 13 regions and invert them as well. For completion, we will give all the pictures now in Figures \ref{fig:firstsubregion}, \ref{fig:secondsubregion}, and \ref{fig:thirdsubregion}. It is immediately clear that all the subregions that result after this second iteration are the same as those that resulted after the first iteration. Thus these 13 regions are the only 13 that will ever appear or need to be analyzed. This happens due to a miraculous happenstance: the inversion of the arcs that make up the borders of the $K_{k,\ell}$ regions, when taken mod $1$, map back onto the exact same arcs again. 

In particular, the boundary of $K$ is composed of pieces of the four lines $\Re(z)=\pm 1/2$ and $\Im(z)=\pm 1/2$. These lines, when inverted become the four circles $|z\pm 1|=1$ and $|z\pm \ii|=1$. When reduced modulo $1$, these become the arcs of the circles $|z+\alpha|=1$ for $\alpha\in \{\pm 1, \pm \ii, \pm 1\pm \ii\}$. The circles $|z+\alpha|=1$ for $\alpha\in \{\pm 1, \pm \ii\}$ when inverted go back to the lines $\Re(z)=\pm 1/2$ or $\Im(z)=\pm 1/2$. The circles $|z+\alpha|=1$ for $\alpha\in \{\pm 1\pm \ii\}$ when inverted are all mapped back onto themselves. This is easily checked by confirming that both they and their inversions must intersect the unit circle at the origin twice (at two of the points $\pm 1, \pm \ii$); and if we check a third point, say the point $1+2\ii$ on the circle $|z-(1+\ii)|=1$, then when inverted this becomes the point $.2-.4\ii$ which lies on the circle $|z-(1-\ii)|=1$. Thus we will only ever see segments of the four lines $\Re(z)=\pm 1/2$, $\Im(z)=\pm 1/2$ and arcs of the $8$ circles  $|z+\alpha|=1$ for $\alpha\in \{\pm 1, \pm \ii, \pm 1\pm \ii\}$.

We may now consider these pictures in the following way: if $z_{i-1}$ belongs to the subregion on the left (in Figures \ref{fig:Kinitial}, \ref{fig:firstsubregion}, \ref{fig:secondsubregion}, or \ref{fig:thirdsubregion}) then the allowable digits $a_{i}$ are the digits corresponding to tiles which intersect the region on the right of the same figure. Moreover, the possible region for $z_{i}$ (given a particular region for $z_{i-1}$ and choice of $a_{i}$) is the corresponding tile shifted back to be centered at the origin.

It is clear from a visual inspection that for most $a_i$ in $G$, the corresponding region for $z_i$ is always the same. For example, the digit $a_i=2$ is allowable in the pictures given by Figures   \ref{fig:Kinv}, \ref{fig:D3inv}, \ref{fig:D4inv}, \ref{fig:D5inv}, \ref{fig:D12inv}, and \ref{fig:D13inv}, and all of Figure \ref{fig:secondsubregion}. In every case, the resulting possible region for $z_i$ is \ref{fig:D3}. 

We see that the only cases where $a_i$ can have multiple corresponding regions for $z_i$ are when $a_i = \pm 2 \pm \ii, \pm 1 \pm 2\ii,\pm2\pm2\ii$, and in each case there are only two possibilities for the corresponding region for $z_i$. We therefore will add to our collection of possible digits a selection of ``marked" digits. For each of these 12 digits, for example $2+\ii $, we will also include the marked digit $(2+\ii )'$. One should consider these marked digits in a similar way as one considers colored partitions. Both $2+\ii $ and $(2+\ii )'$ have the same numerical value for computing the continued fraction but exist to help separate the two possible regions for the next $z_i$. In particular, we will let the unmarked digit, such as $2+\ii $, be such that if $a_i =2+\ii $, then the allowable region for $z_i$ is the same one that would appear had the allowable region for $z_{i-1}$ been Figure \ref{fig:Kinv}, in this case, the region would be given by that of Figure \ref{fig:D9}. For the marked digit, such as $(2+\ii )'$, we obtain the other allowable region for the next $z_i$: in the example of $a_i=(2+\ii )'$, the allowable region for $z_i$ will be Figure \ref{fig:D3}, which appears as a result of the allowable region for $z_{i-1}$ being Figure \ref{fig:D8inv}.

We let $G'=G\cup \{(\pm 2\pm \ii)', (\pm 1\pm 2\ii)', (\pm2\pm2\ii)'\}$. With the definitions of the marked digits, we see that each $a_i\in G'$ corresponds to a uniquely defined allowable region for $z_i$, and thus to a uniquely defined set of possible $a_{i+1}$'s. We can thus describe all allowable strings of digits in the Hurwitz continued fraction expansion by a simple one-step process: that is, the only thing that is needed to understand whether a digit can occur in a given position is to know what the previous digit was. This is the content of Theorem \ref{thm:admissible}.

As such, we have the table given in Table \ref{table:allowabledigits}.
\begin{table}
\begin{tabular}{ccc}
Value of  $a_i$ & Corresponding region for $ z_i$ & Possible values of $a_{i+1}$\\
$2$, $(2+\ii)'$ & Figure \ref{fig:D3} & $\{x+y \ii\in G: |x|\ge 0\}$\\
$2\ii$, $(1+2\ii)'$ & Figure \ref{fig:D5} & $\{x+y \ii\in G: |y|\le 0\}$\\
$1+\ii$ & Figure \ref{fig:D13} & $\{x+y \ii\in G: |x|\ge 0, |y|\le 0\}$\\
$2+\ii$, $1+2\ii$, or $(2+2\ii)'$  & Figure \ref{fig:D9} & $G \setminus \{-1+\ii,-1+2\ii,-2+\ii,-2+2\ii\} $\\
 & & $\cup \{(-1+2\ii)',(-2+\ii)',(-2+2\ii)'\}$
\end{tabular}
\caption{Admissible digit successors}
\label{table:allowabledigits}
\end{table}
We obtain the necessarily relations for the remaining digits in the following way: if we take the negative and/or conjugate on the left-hand column, we also take negatives and/or conjugates in the right-hand column. For all remaining possible values of $a_i$, the corresponding region for $z_i$ is given by Figure \ref{fig:K} and the possible values of $a_{i+1}$ are everything in $G$.

To create an allowable string for the Hurwitz expansion, one follows the above rules and then, erases any marks from marked digits.

\section{Our method}\label{sec:method}

We will now detail the program we used to achieve our calculations.

Let us emphasize to begin with that our method is really based on the following two-step process. First, we generate an array of boolean ``pixels" that represent a rough plot of the points contained in a given $V_{k,\ell}$. (We will throughout this section consider $V_{1,1}$ for simplicity.) Second, to approximate the Taylor series of $h(z)$, we approximate the integral in \eqref{eq:hderivative} by a sum over the points in our array of pixels. Again, for simplicity, we will approximate the Taylor series at the point $-.5-.5\ii$. In other words, we are calculating $h_{m,n}=h_{m,n}^{1,1}(-.5,-.5)$.

The choice of $-.5-.5\ii$ was initially chosen since it was a rational point in $V_{1,1}$ with small denominator (a ``simple" point in some sense), but as it turned out, it was a very good choice to make and greatly simplified the resulting calculations, see Conjecture \ref{conj:odd}.

Let us expand this and work through our method in more thorough detail.

\textbf{Step 1: Initialization}

First, we need a point whose orbit we believe to be dense. We will use the point $(z_0,w_0)=((\log 4-1)+(\log 7-2)\ii,0)$, which was also used by Hensley in some of his calculations available on his website.

We will choose a positive integer $k$ (not related to the $k$ in $K_{k,\ell}$) that will help represent the degree of approximation we achieve. We will then let $Q=2^k$. Then let \textbf{PixelArray} denote a $2Q\times 2Q$ array of boolean values all set to false to begin with. 

Each point in the array is meant to correspond to a small square in the complex plane. In particular, the $(i,j)$ value is meant to correspond to the square of side-length $1/Q$ centered at 
\[
\frac{i+j\ii}{Q}-(1+\ii)\left( 1+\frac{1}{2Q}\right).
\]
In particular, $(i,j)$ corresponds to the square 
\[
\left[\frac{i-1}{Q}-1,\frac{i}{Q}-1\right]\times \left[\frac{j-1}{Q}-1,\frac{j}{Q}-1\right]
\]
considered by treating $\mathbb{C}$ as $\mathbb{R}^2$ in the usual way.
Note that these squares are all disjoint except for their borders and that their union is the entire square centered at the origin with side-length $2$. In the next step, we will attempt to flip the value at $(i,j)$ from false to true if this square intersects $V_{1,1}$, so if we represented true values of $(i,j)$ as a filled-in black square and false faluse of $(i,j)$ as a filled-in white square, we would essentially have a rough pixel depiction of $V_{1,1}$, hence the name. 

Finally, since we performed our calculations in Mathematica, we found it beneficial to precalculate the function
\begin{equation}\label{eq:Hdef}
H_{m,n}(a,b) = \left. \frac{\partial^{m+n}}{\partial x^m \partial y^n} \frac{1}{((ax-by+1)^2+(ay+bx)^2)^2}\right|_{x=y=-.5}
\end{equation}
as a compiled function. In particular these derivatives quickly would become computationally intensive to recalculate every time we wanted to call the function. To demonstrate this, here are some of the $H_{m,n}(a,b)$ functions:
\begin{align*}
H_{0,0}(a,b)&= \frac{1}{(1-a+a^2/2+b+b^2/2)^2}\\
H_{0,1}(a,b)&= \frac{2(-a^2-2b-b^2)}{(1-a+a^2/2+b+b^2/2)^3}\\
H_{0,2}(a,b)&= \frac{2(-a^2+a^3+a^4+5a^2b+5b^2+ab^2+2a^2b^2+5b^3+b^4)}{(1-a+a^2/2+b+b^2/2)^4},
\end{align*}
and so on.

In fact, to save on time, we generated a matrix $\mathcal{H}(a,b)=(H_{m,n}(a,b))_{m,n=1}^L$ for some choice of $L$.

For studying $V_{1,1}$ around $-.5-.5\mathbf{i}$, we could use the natural symmetries of $V_{1,1}$ to only calculate the lower-triangular matrix of $\mathcal{H}(a,b)$ isntead of the full matrix.  We will ultimately sum $H_{m,n}$ over values of $(a,b)$ corresponding to the center of true pixels in \textbf{PixelArray}. Since $V_{1,1}$ is symmetric around $a=-b$ (seeing $w=a+b\mathbf{i}$), we could replace $a$ with $-b$ (and vice-versa), in the sum. However, if in \eqref{eq:Hdef} we replace $a$ with $-b$ and $x$ with $y$, the value of the function is unchanged. Again, replacing $a$ with $-b$ has no effect once we have completed the sum and replacing $x$ with $y$ has the effect of swapping $m$ and $n$, hence why we could, if we desired, look only over the lower-triangular matrix.

This behavior is expected.  Recall that $\mu(\ii A)=\mu(\overline{A})=\mu(A)$ by Theorem \ref{thm:Hensley}. As a consequence $K_{1,1}$ is symmetric along the line $x=y$, and since $-.5-.5\ii$ is on this line, we have that $h_{m,n}=h_{n,m}$.

Note that $V_{2,1}$ is symmetric over $a=-b$ as well, but $V_{3,1}$ is instead symmetric over $b=0$. See Figure \ref{fig:KVdecomp}.

\textbf{Step 2: Populating the pixel array}

We will illustrate the simple method to do this in this section, and later on will give a more complex (but in some cases, time-saving) method later on.

Let $(z,w)=(z_0,w_0)$. Perform the following operation $100\times Q^2$ times: Let $(z,w):=\hat{T}(z,w)$ in the usual assignment sense, and then if $z\in K_{1,1}$, let the $(i,j)$ coordinate of \textbf{PixelArray} be true, where 
\[
i=\lfloor \Re(w)Q+1\rfloor \qquad \text{ and } \qquad j=\lfloor \Im(w)Q+1\rfloor
\]
In particular, $w$ will belong to the square that coordesponds to the coordinates $(i,j)$.

We note that $100*Q^2$ here is somewhat arbitrary. It was merely a value we used that seemed to fill the array rather well.

\textbf{Step 3: (Optional) Fill the holes in the array}

(This step is entirely optional and if not used may just be skipped over completely.)

As was shown in \cite{EINN}, the region $V_{1,1}$ is simply connected and should not contain any ``holes." However, our pixel array, as it is generated in an ad hoc and somewhat probabilitistic fashion, might contain holes, by which we mean pixels $(i,j)$ that are false despite the entire corresponding square belonging to $V_{1,1}$. 

There are a number of methods one could employ to fill in these holes. We illustrate a few here.

\begin{enumerate}
\item The simplest method would be to make use of the natural symmetry inherent in $V_{1,1}$. In particular, it is symmetric around the line $a=-b$. Thus, one could create a new array of pixels where $(i,j)$ is true if and only if either $(i,j)$ or $(2Q+1-j,2Q+1-i)$ is true in the original array of pixels.
\item One could alternately tweak the original array of pixels by making $(i,j)$ be true if its neighbors, $(i+1,j)$, $(i-1,j)$, $(i,j+1)$, and $(i,j-1)$, are all true. 
\item By far the most effective but also most time- and memory-intensive method is to use a flood-fill algorithm.

In particular, one creates a new $2Q\times 2Q$ pixel array, let us call it \textbf{NewPixelArray}, that is false in every coordinate. Then one turns the $(1,1)$ coordinate of \textbf{NewPixelArray} true and performs a flood-fill algorithm: we look above, below, to the left and right of every pixel in \textbf{NewPixelArray} that we have turned true and also turn those true unless the corresponding pixel in \textbf{PixelArray} is true. When completed, this will give us essentially a photo-negative of the desired array. We can then let \textbf{PixelArray} be \textbf{NewPixelArray} after applying negation in every coordinate.
\end{enumerate}

Our preferred method was to mix the flood-fill algorithm with facts about symmetry. Flood-fill algorithms are inherently memory-intensive, so we implemented the following three optimizations: we first used a scanline variant which saves considerably on memory, we made use of the symmetry of $V_{1,1}$ to realize that we only needed to calculate \textbf{NewPixelArray} for $i\le 2Q+1-j$ (i.e., the lower-triangular portion), and then we made use of the symmetry a second time to compare $(i,j)$ of \textbf{NewPixelArray} to both $(i,j)$ and $(2Q+1-j,2Q+1-i)$ of \textbf{PixelArray}, choosing to turn $(i,j)$ true in \textbf{NewPixelArray} unless \emph{either} of the pixels in \textbf{PixelArray} is true.

Finally, we note that regardless of which process we use, we might accidentally fill in holes that we were not supposed to, due to the fractal nature of the boundary; however, in various experiments we felt that the benefit from doing this outweighed the downsides.

\textbf{Step 4: Approximate the invariant measure}

Now, we wish to approximate
\[
h_{m,n}=\left.\frac{\partial^{m+n}}{\partial^m x\partial^n y} h(x,y)\right|_{x=y=-.5}.
\]

We will do this quite simply, by taking smoothed sum over $\mathcal{H}(a,b)$ over $(a,b)$ corresponding to the centers of true pixels $(i,j)$. In particular, we approximate it by
\begin{equation}\label{eq:Hmnsum}
\sum_{\substack{(i,j)\in \textbf{PixelArray}\\(i,j)\text{ true}}} \frac{F(i,j)}{Q^2}\cdot H_{m,n}(a_i,b_j),
\end{equation}
where 
\[
a_i=\frac{i}{Q}-\frac{1}{2Q}-1, \qquad b_j=\frac{j}{Q}-\frac{1}{2Q}-1,
\]
and $F_{i,j}$ is the number of the pixels $(i+\alpha,j+\beta)$, $\alpha,\beta\in \{-1,0,1\}$ that are true divided by $9$. Loosely speaking, $F_{i,j}$ can be seen as a generalization of the trapezoidal rule of numerical integration. It's presence is not necessary for the method to work, but we found it promoted faster convergence. 

To save time, we could replace $H_{m,n}$ in \eqref{eq:Hmnsum} with the matrix $\mathcal{H}$.

\subsection{A variant on step 2}

In step 2, as described above, even after running $100*Q^2$ iterations, we were often left with a large number of holes that needed to be filled, and even then the boundary was not as precise as we would like. This is for two reasons. 

First, $z$ only visits $K_{1,1}$ a small portion of the time, and hence the information from lots of iterates is not being used. Mathematica calculations suggest that $\mu(K_{1,1})=0.066$, so that $z$ visits $K_{1,1}$ roughly $1/15$ of the time.

Second, we spend a long time filling in the middle of $V_{1,1}$ when this information is completely ignored when we use the flood-fill algorithm.

\begin{figure}[h!]
\centering
\includegraphics[width=15cm]{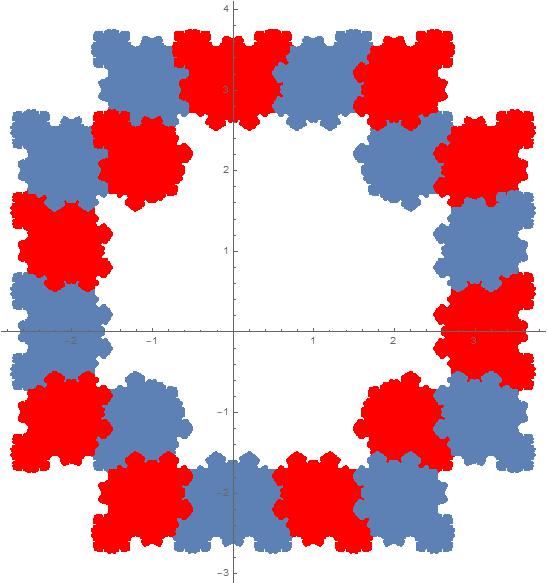}
\caption{The boundary of $V_{1,1}^{-1}$, showing the translates of the $V_{k,\ell}$'s that compose it.}
\label{fig:V11inverse}

\end{figure}

We can mitigate both of these problems using a slightly different method. Let $V_{1,1}^{-1}$ denote the set of points $1/w$ for $w\in V_{1,1}$. We know that $\hat{T}(z,w)=(Tz,1/(w+\alpha))$ for some integer $\alpha\in \mathbb{Z}[\mathbbm{i}]$ depending only on $z$. So let us abuse notation for a moment. If $z\in K_{k,\ell}$, then $\hat{T}(z,V_{k,\ell})=(Tz,1/(V_{k,\ell}+\alpha))$. So if $Tz\in K_{1,1}$, then $(V_{k,\ell}+\alpha)^{-1}\subset V_{1,1}$. In particular, this tells us that $V_{1,1}^{-1}$ is a union of translates of the $V_{k,\ell}$'s. In fact, from \cite{EINN}, it can be shown that this is a disjoint union,  and the boundary of $V_{1,1}^{-1}$ intersects finitely many of these translates. This can be seen directly in Figure \ref{fig:V11inverse}.

So what we do is, like in step 2 above, perform the assignment $(z,w):=\hat{T}(z,w)$ a large number of times (we found $3*Q^2$ was sufficient), and now at every iteration, we calculate the $k,\ell$ such that $z\in K_{k,\ell}$, find all the $\alpha$'s such that the translate $V_{k,\ell}+\alpha$ is a part of the boundary of $ V_{1,1}^{-1}$, and then turn each of the $(i,j)$'s in \textbf{PixelArray} true that correspond to the points $1/(w+\alpha)$ for each of these $\alpha$'s. The values of the $\alpha$'s for each $k$ and $\ell$ can be precalculated to speed up this process: this can be accomplished by using the admissible digits sequences studied previously, although for practical purposes, they can be read directly from pictures such as Figure \ref{fig:V11inverse}. A table of the corresponding $\alpha$'s for $V_{1,1}$ can be seen in Table \ref{tab:V11alpha}.

\begin{table}[h!]

\caption{For a given, $K_{k,\ell}$, the $\alpha\in\mathbb{Z}[\ii]$ such that $V_{k,\ell}+\alpha$ is part of the boundary of $V_{1,1}^{-1}$.}\label{tab:V11alpha}

\begin{tabular}{|c|c|}
\hline
$K_{k,\ell}$ & $\alpha$\\
\hline\hline
$K_{1,1}$ & $-1+2\ii$\\
\hline
$K_{2,1}$ & $-2+\ii$, $-2+2\ii$, $-1+3\ii$\\
\hline
$K_{3,1}$ & $-2$\\
\hline
$K_{1,2}$ & $2+2\ii$\\
\hline
$K_{2,2}$ & $1+2\ii$, $2+2\ii$, $3+2\ii$, $3+\ii$\\
\hline
$K_{3,2}$ & $3\ii$\\
\hline
$K_{1,3}$ & $2-\ii$\\
\hline
$K_{2,3}$ & $1-2\ii$, $2-2\ii$, $3-\ii$\\
\hline
$K_{3,3}$ & $3$\\
\hline
$K_{1,4}$ & $-1-\ii$\\
\hline
$K_{2,4}$ & $-1-2\ii$, $-2-\ii$\\
\hline
$K_{3,4}$ & $-2\ii$\\
\hline
\end{tabular}
\end{table}

For example, by examining Figure \ref{fig:V11inverse}, we can see that if $z\in K_{1,1}$, then since the only copy of $V_{1,1}$ that appears in this image has been translated by $-1+2\mathbbm{i}$, we  have that $(z-1+2\mathbbm{i})^{-1}\in V_{1,1}$. If $z\in K_{2,1}$, then since the copies of $V_{2,1}$ that appear in this image have been translated by $-2+\ii$, $-2+2\ii$, and $-1+3\ii$, we have that $(z-2+\ii)^{-1},(z-2+2\ii)^{-1},(z-1+3\ii)^{-1}\in V_{1,1}$. And so on.

We can even optimize this further, since if we know that $w\in V_{k,\ell}$ for a given $k,\ell$, then $(-\mathbbm{i})^j w \in \cdot V_{k,\ell+j\mod{4}}$ for $j=1,2,3$.

\section{Results}\label{sec:results}

We performed several iterations of the methods listed above. We ran calculations from $k=7$ up to $k=13$. In our implementation, for $V_{1,1}$ centered at $x+y \ii=-.5-.5\ii$, we were able to calculate approximations to $h_{m,n}=h_{m,n}^{1,1}(-.5,-.5)$ for $0\le m,n\le 8$ with $k=7$ in three seconds and using $7.5$ megabytes of memory. The same calculations run for $k=13$ took $10840$ seconds, around three and a half hours, and just over $22$ gigabytes of memory. Increasing $k$ by $1$ typically increased both the memory required and the running time by a factor of a little more than $4$.

Clearly, the biggest obstacle to extending these results further is the amount of fast memory used and the primary contributor to the memory used is the flood-fill algorithm. More precise calculations will likely require slower but less memory-intensive methods, likely dropping any form of filling algorithm and running more iterations in Step 1 of the procedure instead. 

Table \ref{tab:h00h01} shows the $h_{m,n}=h_{m,n}^{1,1}(-.5,-.5)$ values we calculated for  $k$ ranging from $k=7$ to $k=13$. What we note from these is that they show a very regular convergence that appears to be exponential in $k$---that is, that the behavior of the approximation to $h_{m,n}$ in $k$ is like $a+b\cdot c^k$ for appropriate variables $a,b,c$ with $|c|<1$. The value of $a$ that we calculate should be a better approximation to the true value of $h_{m,n}$. We used a least-squares approximation (via Mathematica's \textbf{FindFit} function) to estimate these $a,b,c$ and compiled the results in Table \ref{tab:V11abc}. In Table \ref{tab:V11abc}, we only calculate $a,b,c$ for those $h_{m,n}$ where $m,n$ are both even, because the others did not display very good behavior for the methods used by \textbf{FindFit} to converge.

In fact, by examining Table \ref{tab:h00h01} more closely, we see that $h_{m,n}$ is very, very small whenever $m$ or $n$ is odd. This leads us to the following conjecture.

\begin{table}
\centering
\caption{Approximations of $h_{m,n}^{1,1}(-.5,-.5)$  with various different $k$'s}\label{tab:h00h01}

\begin{tabular}{|c||c|c|c|c|c|c|c|c|}
\hline
$k$ & $7$ & $8$ & $9$ & $10$ & $11$ & $12$ & $13$  \\
\hline

$h_{0,0}$

& $0.76153$ 

& $0.74149$ 

& $0.73038$ 

& $0.72341$ 

& $0.71963$ 

& $0.71732$ 

& $0.71608$ 

\\

\hline

$h_{0,1}$

& $-0.01852$ 

& $-0.00992$ 

& $-0.00695$ 

& $-0.00376$ 

& $-0.00222$ 

& $-0.00120$ 

& $-0.00071$ 

\\

\hline

$h_{0,2}$

& $0.38551$ 

& $0.36624$ 

& $0.35565$ 

& $0.34916$ 

& $0.34562$ 

& $0.34346$ 

& $0.34231$ 

\\

\hline

$h_{0,3}$

& $-0.01156$ 

& $-0.00595$ 

& $-0.00413$ 

& $-0.00223$ 

& $-0.00133$ 

& $-0.00071$ 

& $-0.00043$ 

\\

\hline

$h_{0,4}$

& $0.10541$ 

& $0.09730$ 

& $0.09309$ 

& $0.09053$ 

& $0.08917$ 

& $0.08836$ 

& $0.08793$ 

\\

\hline

$h_{0,5}$

& $-0.00397$ 

& $-0.00199$ 

& $-0.00136$ 

& $-0.00071$ 

& $-0.00043$ 

& $-0.00022$ 

& $-0.00014$ 

\\

\hline

$h_{0,6}$

& $0.02365$ 

& $0.02098$ 

& $0.01969$ 

& $0.01891$ 

& $0.01851$ 

& $0.01828$ 

& $0.01816$ 

\\

\hline

$h_{0,7}$

& $-0.00161$ 

& $-0.00074$ 

& $-0.00050$ 

& $-0.00025$ 

& $-0.00015$ 

& $-0.00008$ 

& $-0.00005$ 

\\

\hline

$h_{0,8}$

& $0.00963$ 

& $0.00867$ 

& $0.00814$ 

& $0.00783$ 

& $0.00766$ 

& $0.00757$ 

& $0.00751$ 

\\

\hline

$h_{1,1}$

& $-0.00716$ 

& $-0.00433$ 

& $-0.00037$ 

& $-0.00064$ 

& $-0.00045$ 

& $-0.00022$ 

& $-0.00016$ 

\\

\hline

$h_{1,2}$

& $-0.01946$ 

& $-0.01107$ 

& $-0.00757$ 

& $-0.00418$ 

& $-0.00237$ 

& $-0.00133$ 

& $-0.00075$ 

\\

\hline

$h_{1,3}$

& $-0.00458$ 

& $-0.00294$ 

& $-0.00023$ 

& $-0.00042$ 

& $-0.00032$ 

& $-0.00015$ 

& $-0.00011$ 

\\

\hline

$h_{1,4}$

& $-0.00774$ 

& $-0.00427$ 

& $-0.00283$ 

& $-0.00144$ 

& $-0.00077$ 

& $-0.00042$ 

& $-0.00024$ 

\\

\hline

$h_{1,5}$

& $-0.00188$ 

& $-0.00018$ 

& $0.00019$ 

& $0.00002$ 

& $0.00002$ 

& $-0.00002$ 

& $-5.87870*10^{-6}$ 

\\

\hline

$h_{1,6}$

& $0.00028$ 

& $0.00018$ 

& $0.00019$ 

& $0.00016$ 

& $0.00011$ 

& $0.00007$ 

& $0.00004$ 

\\

\hline

$h_{1,7}$

& $-0.00021$ 

& $0.00029$ 

& $0.00015$ 

& $0.00007$ 

& $0.00005$ 

& $9.30540*10^{-6}$ 

& $0.00001$ 

\\

\hline

$h_{1,8}$

& $0.00185$ 

& $0.00039$ 

& $0.00021$ 

& $0.00009$ 

& $0.00004$ 

& $0.00003$ 

& $0.00001$ 

\\

\hline

$h_{2,2}$

& $0.58500$ 

& $0.54861$ 

& $0.52742$ 

& $0.51478$ 

& $0.50756$ 

& $0.50306$ 

& $0.50069$ 

\\

\hline

$h_{2,3}$

& $-0.02677$ 

& $-0.01484$ 

& $-0.01021$ 

& $-0.00586$ 

& $-0.00339$ 

& $-0.00192$ 

& $-0.00110$ 

\\

\hline

$h_{2,4}$

& $0.32479$ 

& $0.29977$ 

& $0.28537$ 

& $0.27685$ 

& $0.27201$ 

& $0.26901$ 

& $0.26744$ 

\\

\hline

$h_{2,5}$

& $-0.00519$ 

& $-0.00404$ 

& $-0.00277$ 

& $-0.00178$ 

& $-0.00109$ 

& $-0.00060$ 

& $-0.00036$ 

\\

\hline

$h_{2,6}$

& $-0.02441$ 

& $-0.03010$ 

& $-0.03116$ 

& $-0.03158$ 

& $-0.03175$ 

& $-0.03181$ 

& $-0.03184$ 

\\

\hline

$h_{2,7}$

& $-0.00208$ 

& $0.00081$ 

& $0.00088$ 

& $0.00050$ 

& $0.00028$ 

& $0.00019$ 

& $0.00011$ 

\\

\hline

$h_{2,8}$

& $-0.04454$ 

& $-0.04429$ 

& $-0.04288$ 

& $-0.04190$ 

& $-0.04132$ 

& $-0.04094$ 

& $-0.04074$ 

\\

\hline

$h_{3,3}$

& $-0.00597$ 

& $-0.00743$ 

& $-0.00131$ 

& $-0.00124$ 

& $-0.00098$ 

& $-0.00034$ 

& $-0.00027$ 

\\

\hline

$h_{3,4}$

& $-0.03571$ 

& $-0.01854$ 

& $-0.01263$ 

& $-0.00684$ 

& $-0.00383$ 

& $-0.00220$ 

& $-0.00123$ 

\\

\hline

$h_{3,5}$

& $-0.00513$ 

& $-0.00433$ 

& $-0.00071$ 

& $-0.00071$ 

& $-0.00057$ 

& $-0.00020$ 

& $-0.00016$ 

\\

\hline

$h_{3,6}$

& $-0.02558$ 

& $-0.00677$ 

& $-0.00371$ 

& $-0.00150$ 

& $-0.00070$ 

& $-0.00038$ 

& $-0.00017$ 

\\

\hline

$h_{3,7}$

& $0.01199$ 

& $0.00259$ 

& $0.00076$ 

& $0.00041$ 

& $0.00030$ 

& $0.00008$ 

& $0.00005$ 

\\

\hline

$h_{3,8}$

& $-0.00455$ 

& $0.00190$ 

& $0.00203$ 

& $0.00114$ 

& $0.00067$ 

& $0.00043$ 

& $0.00025$ 

\\

\hline

$h_{4,4}$

& $0.68944$ 

& $0.64689$ 

& $0.61713$ 

& $0.59915$ 

& $0.58883$ 

& $0.58232$ 

& $0.57890$ 

\\

\hline

$h_{4,5}$

& $-0.01286$ 

& $-0.01699$ 

& $-0.01276$ 

& $-0.00766$ 

& $-0.00451$ 

& $-0.00262$ 

& $-0.00153$ 

\\

\hline

$h_{4,6}$

& $0.31911$ 

& $0.29305$ 

& $0.27634$ 

& $0.26633$ 

& $0.26067$ 

& $0.25713$ 

& $0.25527$ 

\\

\hline

$h_{4,7}$

& $0.00691$ 

& $-0.00096$ 

& $-0.00079$ 

& $-0.00130$ 

& $-0.00095$ 

& $-0.00051$ 

& $-0.00033$ 

\\

\hline

$h_{4,8}$

& $-0.09596$ 

& $-0.10365$ 

& $-0.10263$ 

& $-0.10041$ 

& $-0.09941$ 

& $-0.09870$ 

& $-0.09825$ 

\\

\hline

$h_{5,5}$

& $-0.03693$ 

& $-0.01355$ 

& $-0.00289$ 

& $-0.00220$ 

& $-0.00171$ 

& $-0.00054$ 

& $-0.00039$ 

\\

\hline

$h_{5,6}$

& $-0.03640$ 

& $-0.02237$ 

& $-0.01594$ 

& $-0.00806$ 

& $-0.00441$ 

& $-0.00261$ 

& $-0.00146$ 

\\

\hline

$h_{5,7}$

& $-0.01519$ 

& $-0.00587$ 

& $-0.00114$ 

& $-0.00098$ 

& $-0.00077$ 

& $-0.00025$ 

& $-0.00018$ 

\\

\hline

$h_{5,8}$

& $-0.04348$ 

& $-0.00868$ 

& $-0.00369$ 

& $-0.00110$ 

& $-0.00031$ 

& $-0.00006$ 

& $1.99160*10^{-6}$ 

\\

\hline

$h_{6,6}$

& $0.65399$ 

& $0.62170$ 

& $0.59129$ 

& $0.57018$ 

& $0.55883$ 

& $0.55159$ 

& $0.54765$ 

\\

\hline

$h_{6,7}$

& $0.01380$ 

& $-0.01494$ 

& $-0.01345$ 

& $-0.00802$ 

& $-0.00483$ 

& $-0.00298$ 

& $-0.00173$ 

\\

\hline

$h_{6,8}$

& $0.26594$ 

& $0.24508$ 

& $0.22976$ 

& $0.21971$ 

& $0.21429$ 

& $0.21087$ 

& $0.20903$ 

\\

\hline

$h_{7,7}$

& $-0.07041$ 

& $-0.01733$ 

& $-0.00461$ 

& $-0.00327$ 

& $-0.00261$ 

& $-0.00071$ 

& $-0.00053$ 

\\

\hline

$h_{7,8}$

& $-0.03244$ 

& $-0.02434$ 

& $-0.01788$ 

& $-0.00847$ 

& $-0.00448$ 

& $-0.00274$ 

& $-0.00154$ 

\\

\hline

$h_{8,8}$

& $0.55655$ 

& $0.54007$ 

& $0.51596$ 

& $0.49254$ 

& $0.48116$ 

& $0.47381$ 

& $0.46970$ 

\\

\hline
\end{tabular}
\end{table}

\begin{table}
\centering
\caption{Approximations to $h_{m,n}^{1,1}(-.5,-.5)$ as in Table \ref{tab:h00h01} of the form $a+b\cdot c^k$.}
\label{tab:V11abc}

\begin{tabular}{|c|c|c|c|}
\hline

$h_{m,n}$ & $a$ & $b$ & $c$ \\

\hline\hline

$h_{0,0}$ & $0.7149$ & $2.6050$ & $0.5633$\\

\hline

$h_{0,2}$ & $0.3411$ & $2.4500$ & $0.5641$\\

\hline

$h_{0,4}$ & $0.0875$ & $1.1610$ & $0.5510$\\

\hline

$h_{0,6}$ & $0.0180$ & $0.4290$ & $0.5380$\\

\hline

$h_{0,8}$ & $0.0075$ & $0.1423$ & $0.5505$\\

\hline

$h_{2,2}$ & $0.4974$ & $3.8590$ & $0.5825$\\

\hline

$h_{2,4}$ & $0.2652$ & $2.7210$ & $0.5793$\\

\hline

$h_{2,6}$ & $-0.0321$ & $5.5710$ & $0.3877$\\

\hline

$h_{2,8}$ & $-0.0409$ & $-0.1466$ & $0.6040$\\

\hline

$h_{4,4}$ & $0.5754$ & $4.8780$ & $0.5864$\\

\hline

$h_{4,6}$ & $0.2532$ & $3.0290$ & $0.5798$\\

\hline

$h_{4,8}$ & $-0.1003$ & $0.1056$ & $0.5410$\\

\hline

$h_{6,6}$ & $0.5484$ & $5.7890$ & $0.5695$\\

\hline

$h_{6,8}$ & $0.2085$ & $3.1280$ & $0.5679$\\

\hline

$h_{8,8}$ & $0.4762$ & $6.8400$ & $0.5399$\\ 
\hline
\end{tabular}
\end{table}

\begin{table}
\centering
\caption{Approximations of $h_{m,n}^{2,1}(0,0)$ with various different $k$'s}
\label{tab:V21table}

\begin{tabular}{|c||c|c|c|c|c|c|c|}
\hline
$k$ & $7$ & $8$ & $9$ & $10$ & $11$ & $12$ & $13$  \\
\hline

$h_{0,0}$ & $1.05000$ & $1.02930$ & $1.01600$ & $1.00920$ & $1.00530$ & $1.00300$ & $1.00170$ \\

\hline

$h_{0,1}$& $0.52909$ & $0.52283$ & $0.51697$ & $0.51422$ & $0.51277$ & $0.51203$ & $0.51153$ \\

\hline

$h_{0,2}$& $0.91839$ & $0.88992$ & $0.87123$ & $0.86185$ & $0.85657$ & $0.85353$ & $0.85180$ \\

\hline

$h_{0,3}$

& $0.47458$ 

& $0.46152$ 

& $0.45109$ 

& $0.44623$ 

& $0.44350$ 

& $0.44207$ 

& $0.44116$ 

\\

\hline

$h_{0,4}$

& $0.36284$ 

& $0.34377$ 

& $0.33079$ 

& $0.32429$ 

& $0.32057$ 

& $0.31851$ 

& $0.31730$ 

\\

\hline

$h_{0,5}$

& $0.08138$ 

& $0.07348$ 

& $0.06702$ 

& $0.06391$ 

& $0.06202$ 

& $0.06105$ 

& $0.06043$ 

\\

\hline

$h_{0,6}$

& $0.02086$ 

& $0.01297$ 

& $0.00729$ 

& $0.00446$ 

& $0.00269$ 

& $0.00176$ 

& $0.00120$ 

\\

\hline

$h_{1,1}$

& $1.11320$ 

& $1.09760$ 

& $1.08560$ 

& $1.07920$ 

& $1.07600$ 

& $1.07410$ 

& $1.07310$ 

\\

\hline

$h_{1,2}$

& $1.67520$ 

& $1.63920$ 

& $1.61480$ 

& $1.60210$ 

& $1.59590$ 

& $1.59210$ 

& $1.59000$ 

\\

\hline

$h_{1,3}$

& $1.79320$ 

& $1.74940$ 

& $1.71870$ 

& $1.70320$ 

& $1.69530$ 

& $1.69060$ 

& $1.68800$ 

\\

\hline

$h_{1,4}$

& $1.17510$ 

& $1.13490$ 

& $1.10670$ 

& $1.09200$ 

& $1.08440$ 

& $1.07990$ 

& $1.07740$ 

\\

\hline

$h_{1,5}$

& $0.15821$ 

& $0.13868$ 

& $0.12366$ 

& $0.11567$ 

& $0.11098$ 

& $0.10834$ 

& $0.10694$ 

\\

\hline

$h_{1,6}$

& $-0.67308$ 

& $-0.67146$ 

& $-0.67239$ 

& $-0.67308$ 

& $-0.67461$ 

& $-0.67549$ 

& $-0.67576$ 

\\

\hline

$h_{2,2}$

& $4.10900$ 

& $3.99640$ 

& $3.92300$ 

& $3.88730$ 

& $3.86880$ 

& $3.85760$ 

& $3.85140$ 

\\

\hline

$h_{2,3}$

& $5.45500$ 

& $5.30780$ 

& $5.21010$ 

& $5.16300$ 

& $5.13990$ 

& $5.12610$ 

& $5.11820$ 

\\

\hline

$h_{2,4}$

& $5.72110$ 

& $5.53150$ 

& $5.40640$ 

& $5.34570$ 

& $5.31510$ 

& $5.29660$ 

& $5.28630$ 

\\

\hline

$h_{2,5}$

& $3.19890$ 

& $3.05920$ 

& $2.96310$ 

& $2.91830$ 

& $2.89430$ 

& $2.88000$ 

& $2.87220$ 

\\

\hline

$h_{2,6}$

& $-0.17949$ 

& $-0.26498$ 

& $-0.32104$ 

& $-0.34417$ 

& $-0.36139$ 

& $-0.37206$ 

& $-0.37688$ 

\\

\hline

$h_{3,3}$

& $11.16100$ 

& $10.85000$ 

& $10.64200$ 

& $10.54200$ 

& $10.49400$ 

& $10.46400$ 

& $10.44700$ 

\\

\hline

$h_{3,4}$

& $15.68100$ 

& $15.19900$ 

& $14.88000$ 

& $14.72400$ 

& $14.65000$ 

& $14.60600$ 

& $14.58000$ 

\\

\hline

$h_{3,5}$

& $15.24300$ 

& $14.70500$ 

& $14.34300$ 

& $14.16900$ 

& $14.08400$ 

& $14.03400$ 

& $14.00500$ 

\\

\hline

$h_{3,6}$

& $8.07010$ 

& $7.62960$ 

& $7.33850$ 

& $7.20920$ 

& $7.13640$ 

& $7.09120$ 

& $7.06770$ 

\\

\hline

$h_{4,4}$

& $29.83100$ 

& $28.91500$ 

& $28.30000$ 

& $27.99600$ 

& $27.85800$ 

& $27.77600$ 

& $27.72700$ 

\\

\hline

$h_{4,5}$

& $41.33600$ 

& $39.98400$ 

& $39.06200$ 

& $38.60000$ 

& $38.39500$ 

& $38.27600$ 

& $38.20300$ 

\\

\hline

$h_{4,6}$

& $39.34600$ 

& $37.84100$ 

& $36.81100$ 

& $36.30700$ 

& $36.07300$ 

& $35.93600$ 

& $35.85400$ 

\\

\hline

$h_{5,5}$

& $78.91400$ 

& $76.32900$ 

& $74.53800$ 

& $73.59900$ 

& $73.20700$ 

& $72.99000$ 

& $72.84900$ 

\\

\hline

$h_{5,6}$

& $109.93000$ 

& $106.20000$ 

& $103.50000$ 

& $102.06000$ 

& $101.47000$ 

& $101.15000$ 

& $100.93000$ 

\\

\hline

$h_{6,6}$

& $216.57000$ 

& $209.70000$ 

& $204.22000$ 

& $201.16000$ 

& $199.96000$ 

& $199.35000$ 

& $198.92000$ 

\\
\hline
\end{tabular}
\end{table}

\begin{table}
\centering
\caption{Approximations of $h_{m,n}^{2,1}(-.5,-.5)$  with various different $k$'s}
\label{tab:V21table2}

\begin{tabular}{|c||c|c|c|c|c|c|c|}
\hline
$k$ & $7$ & $8$ & $9$ & $10$ & $11$ & $12$ & $13$  \\
\hline

$h_{0,0}$

& $1.02990$ 

& $1.00040$ 

& $0.98326$ 

& $0.97423$ 

& $0.96892$ 

& $0.96570$ 

& $0.96398$ 

\\

\hline

$h_{0,1}$

& $-0.47009$ 

& $-0.44193$ 

& $-0.42609$ 

& $-0.41790$ 

& $-0.41280$ 

& $-0.40966$ 

& $-0.40802$ 

\\

\hline

$h_{0,2}$

& $0.86292$ 

& $0.81620$ 

& $0.78943$ 

& $0.77567$ 

& $0.76728$ 

& $0.76222$ 

& $0.75952$ 

\\

\hline

$h_{0,3}$

& $-0.50773$ 

& $-0.46589$ 

& $-0.44346$ 

& $-0.43147$ 

& $-0.42418$ 

& $-0.41985$ 

& $-0.41758$ 

\\

\hline

$h_{0,4}$

& $0.65407$ 

& $0.60225$ 

& $0.57285$ 

& $0.55735$ 

& $0.54796$ 

& $0.54257$ 

& $0.53969$ 

\\

\hline

$h_{0,5}$

& $-0.56716$ 

& $-0.51938$ 

& $-0.49097$ 

& $-0.47590$ 

& $-0.46668$ 

& $-0.46139$ 

& $-0.45862$ 

\\

\hline

$h_{0,6}$

& $0.52681$ 

& $0.47936$ 

& $0.45076$ 

& $0.43584$ 

& $0.42671$ 

& $0.42154$ 

& $0.41880$ 

\\

\hline

$h_{1,1}$

& $0.90306$ 

& $0.86700$ 

& $0.84281$ 

& $0.83173$ 

& $0.82388$ 

& $0.81909$ 

& $0.81647$ 

\\

\hline

$h_{1,2}$

& $-1.04660$ 

& $-0.98469$ 

& $-0.94699$ 

& $-0.92931$ 

& $-0.91688$ 

& $-0.90900$ 

& $-0.90481$ 

\\

\hline

$h_{1,3}$

& $1.17640$ 

& $1.09980$ 

& $1.05580$ 

& $1.03380$ 

& $1.01890$ 

& $1.00970$ 

& $1.00480$ 

\\

\hline

$h_{1,4}$

& $-1.68650$ 

& $-1.56930$ 

& $-1.49670$ 

& $-1.45870$ 

& $-1.43390$ 

& $-1.41900$ 

& $-1.41120$ 

\\

\hline

$h_{1,5}$

& $2.23230$ 

& $2.10270$ 

& $2.00580$ 

& $1.95420$ 

& $1.92010$ 

& $1.89980$ 

& $1.88920$ 

\\

\hline

$h_{1,6}$

& $-2.46230$ 

& $-2.32600$ 

& $-2.21390$ 

& $-2.15380$ 

& $-2.11390$ 

& $-2.09010$ 

& $-2.07770$ 

\\

\hline

$h_{2,2}$

& $1.42810$ 

& $1.34060$ 

& $1.29480$ 

& $1.27360$ 

& $1.25820$ 

& $1.24720$ 

& $1.24150$ 

\\

\hline

$h_{2,3}$

& $-0.55036$ 

& $-0.46931$ 

& $-0.44594$ 

& $-0.43614$ 

& $-0.42732$ 

& $-0.42010$ 

& $-0.41626$ 

\\

\hline

$h_{2,4}$

& $1.58020$ 

& $1.43480$ 

& $1.37120$ 

& $1.33820$ 

& $1.31500$ 

& $1.29910$ 

& $1.29090$ 

\\

\hline

$h_{2,5}$

& $-2.92120$ 

& $-2.79160$ 

& $-2.68540$ 

& $-2.62610$ 

& $-2.58280$ 

& $-2.55440$ 

& $-2.53980$ 

\\

\hline

$h_{2,6}$

& $4.31290$ 

& $4.21970$ 

& $4.07640$ 

& $3.99000$ 

& $3.92570$ 

& $3.88290$ 

& $3.86130$ 

\\

\hline

$h_{3,3}$

& $-1.33430$ 

& $-1.42160$ 

& $-1.37300$ 

& $-1.34430$ 

& $-1.32730$ 

& $-1.31950$ 

& $-1.31580$ 

\\

\hline

$h_{3,4}$

& $1.89940$ 

& $2.10410$ 

& $2.05750$ 

& $2.02980$ 

& $2.01050$ 

& $2.00080$ 

& $1.99630$ 

\\

\hline

$h_{3,5}$

& $-0.17617$ 

& $-0.36837$ 

& $-0.35618$ 

& $-0.35242$ 

& $-0.35211$ 

& $-0.35485$ 

& $-0.35659$ 

\\

\hline

$h_{3,6}$

& $-2.32780$ 

& $-2.27720$ 

& $-2.28260$ 

& $-2.25770$ 

& $-2.23060$ 

& $-2.20410$ 

& $-2.19180$ 

\\

\hline

$h_{4,4}$

& $-4.28560$ 

& $-4.80100$ 

& $-4.73640$ 

& $-4.68480$ 

& $-4.63880$ 

& $-4.60840$ 

& $-4.59460$ 

\\

\hline

$h_{4,5}$

& $4.74050$ 

& $5.40580$ 

& $5.37630$ 

& $5.32910$ 

& $5.28310$ 

& $5.24700$ 

& $5.23160$ 

\\

\hline

$h_{4,6}$

& $-2.29170$ 

& $-2.83350$ 

& $-2.77740$ 

& $-2.74410$ 

& $-2.71530$ 

& $-2.69770$ 

& $-2.69020$ 

\\

\hline

$h_{5,5}$

& $-6.09960$ 

& $-7.15240$ 

& $-7.32090$ 

& $-7.29640$ 

& $-7.25190$ 

& $-7.19250$ 

& $-7.16990$ 

\\

\hline

$h_{5,6}$

& $4.74530$ 

& $5.43940$ 

& $5.76110$ 

& $5.75090$ 

& $5.72840$ 

& $5.67070$ 

& $5.64970$ 

\\

\hline

$h_{6,6}$

& $-4.71690$ 

& $-3.02340$ 

& $-3.86850$ 

& $-3.79360$ 

& $-3.78910$ 

& $-3.68790$ 

& $-3.64450$ 

\\
\hline
\end{tabular}
\end{table}

\begin{table}
\centering
\caption{Approximations of $h_{m,n}^{3,1}(-.5,0)$  with various different $k$'s}
\label{tab:V31table}

\begin{tabular}{|c||c|c|c|c|c|c|c|}
\hline
$k$ & $7$ & $8$ & $9$ & $10$ & $11$ & $12$ & $13$  \\
\hline

$h_{0,0}$

& $1.12250$ 

& $1.09650$ 

& $1.07910$ 

& $1.07040$ 

& $1.06540$ 

& $1.06240$ 

& $1.06080$ 

\\

\hline

$h_{0,2}$

& $0.91900$ 

& $0.88322$ 

& $0.85904$ 

& $0.84794$ 

& $0.84137$ 

& $0.83776$ 

& $0.83561$ 

\\

\hline

$h_{0,4}$

& $0.41936$ 

& $0.39495$ 

& $0.37836$ 

& $0.37117$ 

& $0.36681$ 

& $0.36457$ 

& $0.36313$ 

\\

\hline

$h_{0,6}$

& $0.27205$ 

& $0.25348$ 

& $0.24084$ 

& $0.23545$ 

& $0.23217$ 

& $0.23042$ 

& $0.22933$ 

\\

\hline

$h_{0,8}$

& $0.15823$ 

& $0.14390$ 

& $0.13522$ 

& $0.13149$ 

& $0.12930$ 

& $0.12813$ 

& $0.12741$ 

\\

\hline

$h_{0,10}$

& $0.06755$ 

& $0.05976$ 

& $0.05518$ 

& $0.05324$ 

& $0.05211$ 

& $0.05152$ 

& $0.05115$ 

\\

\hline

$h_{0,12}$

& $0.03038$ 

& $0.02646$ 

& $0.02419$ 

& $0.02322$ 

& $0.02265$ 

& $0.02235$ 

& $0.02216$ 

\\

\hline

$h_{1,0}$

& $-0.04283$ 

& $-0.02533$ 

& $-0.01365$ 

& $-0.00820$ 

& $-0.00478$ 

& $-0.00263$ 

& $-0.00156$ 

\\

\hline

$h_{1,2}$

& $-0.07484$ 

& $-0.04020$ 

& $-0.02038$ 

& $-0.01199$ 

& $-0.00693$ 

& $-0.00380$ 

& $-0.00221$ 

\\

\hline

$h_{1,4}$

& $-0.01383$ 

& $-0.00387$ 

& $-0.00026$ 

& $-0.00036$ 

& $0.00003$ 

& $-0.00011$ 

& $-0.00004$ 

\\

\hline

$h_{1,6}$

& $-0.00015$ 

& $-0.00174$ 

& $0.00038$ 

& $-0.00042$ 

& $-0.00034$ 

& $-0.00024$ 

& $-0.00017$ 

\\

\hline

$h_{1,8}$

& $-0.02442$ 

& $-0.01165$ 

& $-0.00476$ 

& $-0.00284$ 

& $-0.00166$ 

& $-0.00090$ 

& $-0.00054$ 

\\

\hline

$h_{1,10}$

& $-0.01200$ 

& $-0.00534$ 

& $-0.00145$ 

& $-0.00085$ 

& $-0.00041$ 

& $-0.00027$ 

& $-0.00015$ 

\\

\hline

$h_{1,12}$

& $-0.00119$ 

& $-0.00174$ 

& $-0.00033$ 

& $-0.00032$ 

& $-0.00023$ 

& $-0.00016$ 

& $-0.00009$ 

\\

\hline

$h_{2,0}$

& $0.79788$ 

& $0.76122$ 

& $0.73782$ 

& $0.72565$ 

& $0.71870$ 

& $0.71439$ 

& $0.71214$ 

\\

\hline

$h_{2,2}$

& $1.82330$ 

& $1.71800$ 

& $1.65490$ 

& $1.62310$ 

& $1.60550$ 

& $1.59460$ 

& $1.58890$ 

\\

\hline

$h_{2,4}$

& $-0.05398$ 

& $-0.08345$ 

& $-0.09558$ 

& $-0.10245$ 

& $-0.10570$ 

& $-0.10676$ 

& $-0.10762$ 

\\

\hline

$h_{2,6}$

& $-0.50402$ 

& $-0.44407$ 

& $-0.42695$ 

& $-0.41813$ 

& $-0.41414$ 

& $-0.41184$ 

& $-0.41082$ 

\\

\hline

$h_{2,8}$

& $0.38647$ 

& $0.37477$ 

& $0.35858$ 

& $0.35202$ 

& $0.34790$ 

& $0.34480$ 

& $0.34324$ 

\\

\hline

$h_{2,10}$

& $0.28840$ 

& $0.25490$ 

& $0.23269$ 

& $0.22539$ 

& $0.22111$ 

& $0.21889$ 

& $0.21759$ 

\\

\hline

$h_{2,12}$

& $0.13172$ 

& $0.12224$ 

& $0.10761$ 

& $0.10319$ 

& $0.10023$ 

& $0.09874$ 

& $0.09766$ 

\\

\hline

$h_{3,0}$

& $-0.03810$ 

& $-0.02329$ 

& $-0.01275$ 

& $-0.00767$ 

& $-0.00455$ 

& $-0.00247$ 

& $-0.00147$ 

\\

\hline

$h_{3,2}$

& $-0.21693$ 

& $-0.12028$ 

& $-0.06484$ 

& $-0.03720$ 

& $-0.02202$ 

& $-0.01171$ 

& $-0.00688$ 

\\

\hline

$h_{3,4}$

& $-0.12881$ 

& $-0.04430$ 

& $-0.02230$ 

& $-0.00998$ 

& $-0.00444$ 

& $-0.00258$ 

& $-0.00125$ 

\\

\hline

$h_{3,6}$

& $0.21068$ 

& $0.09479$ 

& $0.04629$ 

& $0.02437$ 

& $0.01402$ 

& $0.00716$ 

& $0.00422$ 

\\

\hline

$h_{3,8}$

& $0.02642$ 

& $0.00261$ 

& $-0.00657$ 

& $-0.00540$ 

& $-0.00471$ 

& $-0.00191$ 

& $-0.00140$ 

\\

\hline

$h_{3,10}$

& $-0.12039$ 

& $-0.03183$ 

& $-0.01350$ 

& $-0.00565$ 

& $-0.00190$ 

& $-0.00052$ 

& $-0.00032$ 

\\

\hline

$h_{3,12}$

& $0.00979$ 

& $0.01732$ 

& $0.01188$ 

& $0.00638$ 

& $0.00478$ 

& $0.00227$ 

& $0.00144$ 

\\

\hline

$h_{4,0}$

& $0.35836$ 

& $0.33653$ 

& $0.32200$ 

& $0.31469$ 

& $0.31042$ 

& $0.30773$ 

& $0.30634$ 

\\

\hline

$h_{4,2}$

& $3.06740$ 

& $2.86920$ 

& $2.75120$ 

& $2.69320$ 

& $2.66060$ 

& $2.63950$ 

& $2.62900$ 

\\

\hline

$h_{4,4}$

& $2.25640$ 

& $1.93880$ 

& $1.81690$ 

& $1.75420$ 

& $1.72340$ 

& $1.70700$ 

& $1.69830$ 

\\

\hline

$h_{4,6}$

& $-2.87200$ 

& $-2.68900$ 

& $-2.58070$ 

& $-2.53380$ 

& $-2.50730$ 

& $-2.48780$ 

& $-2.47870$ 

\\

\hline

$h_{4,8}$

& $-1.28790$ 

& $-1.04270$ 

& $-0.92295$ 

& $-0.88697$ 

& $-0.86753$ 

& $-0.86030$ 

& $-0.85562$ 

\\

\hline

$h_{4,10}$

& $-0.00160$ 

& $-0.06368$ 

& $-0.01779$ 

& $-0.00903$ 

& $-0.00027$ 

& $0.00274$ 

& $0.00710$ 

\\

\hline

$h_{4,12}$

& $-1.79360$ 

& $-1.59120$ 

& $-1.46990$ 

& $-1.41090$ 

& $-1.37830$ 

& $-1.35700$ 

& $-1.34480$ 

\\

\hline

$h_{5,0}$

& $-0.01443$ 

& $-0.00960$ 

& $-0.00519$ 

& $-0.00327$ 

& $-0.00196$ 

& $-0.00107$ 

& $-0.00065$ 

\\

\hline

$h_{5,2}$

& $-0.26254$ 

& $-0.15429$ 

& $-0.08308$ 

& $-0.04900$ 

& $-0.02973$ 

& $-0.01564$ 

& $-0.00934$ 

\\

\hline

$h_{5,4}$

& $-0.69504$ 

& $-0.29601$ 

& $-0.14757$ 

& $-0.07613$ 

& $-0.04180$ 

& $-0.02210$ 

& $-0.01259$ 

\\

\hline

$h_{5,6}$

& $0.24273$ 

& $0.14661$ 

& $0.09344$ 

& $0.05626$ 

& $0.03781$ 

& $0.01763$ 

& $0.01134$ 

\\

\hline

$h_{5,8}$

& $0.73041$ 

& $0.21187$ 

& $0.07405$ 

& $0.02813$ 

& $0.00655$ 

& $0.00224$ 

& $0.00083$ 

\\

\hline

$h_{5,10}$

& $-0.35739$ 

& $-0.19373$ 

& $-0.12043$ 

& $-0.06233$ 

& $-0.04104$ 

& $-0.01836$ 

& $-0.01171$ 

\\

\hline

$h_{5,12}$

& $0.25429$ 

& $0.07150$ 

& $0.08992$ 

& $0.05615$ 

& $0.04221$ 

& $0.02215$ 

& $0.01336$ 

\\

\hline

$h_{6,0}$

& $0.05726$ 

& $0.05228$ 

& $0.04809$ 

& $0.04605$ 

& $0.04479$ 

& $0.04398$ 

& $0.04356$ 

\\

\hline

$h_{6,2}$

& $2.18990$ 

& $2.06130$ 

& $1.96470$ 

& $1.91870$ 

& $1.89130$ 

& $1.87250$ 

& $1.86320$ 

\\

\hline

$h_{6,4}$

& $6.68140$ 

& $5.95380$ 

& $5.61530$ 

& $5.45410$ 

& $5.36920$ 

& $5.31670$ 

& $5.29050$ 

\\

\hline

$h_{6,6}$

& $-0.32486$ 

& $-0.65476$ 

& $-0.78060$ 

& $-0.80787$ 

& $-0.81930$ 

& $-0.81160$ 

& $-0.81093$ 

\\

\hline

$h_{6,8}$

& $-2.94040$ 

& $-2.23620$ 

& $-2.12620$ 

& $-2.07260$ 

& $-2.05540$ 

& $-2.04350$ 

& $-2.04500$ 

\\

\hline

$h_{6,10}$

& $7.44640$ 

& $6.57510$ 

& $6.17320$ 

& $5.93830$ 

& $5.81550$ 

& $5.72740$ 

& $5.68180$ 

\\

\hline

$h_{6,12}$

& $2.42410$ 

& $1.48330$ 

& $1.07860$ 

& $1.03600$ 

& $1.00290$ 

& $0.99715$ 

& $0.99556$ 

\\

\hline
\end{tabular}

\end{table}

\begin{conj}\label{conj:odd}
We have $h_{m,n}^{1,1}(-.5,-.5)=0$ whenever $m$ or $n$ is odd.
\end{conj}

First, we emphasize that this conjecture is not a simple consequence of the symmetries of $K_{1,1}$ and $V_{1,1}$. Recall \eqref{eq:hderivative}:
\[
h_{m,n}=h_{m,n}^{1,1}(-.5,-.5)= \int_{a+b\ii\in V_{1,1}}\left.\frac{\partial^{m+n}}{\partial x^m \partial y^n} \frac{1}{((ax-by+1)^2+(ay+bx)^2)^2}\right|_{x+ y\ii =-.5-.5\ii} \ da \ db.
\]
The set $V_{1,1}$ is symmetric under the transformation $a\mapsto -b$, $b\mapsto -a$, so if the integrand at odd $m$ or $n$ is negated by this same transformation, we should get that $h_{m,n}=0$. However, if we examine the derivative inside this equation with, say, $m=1$, $n=0$, then we get
\[
\left.\frac{\partial}{\partial x} \frac{1}{((ax-by+1)^2+(ay+bx)^2)^2} \right|_{x+ y\ii = -.5-.5\ii} =  \frac{-4a+2(a^2+b^2)}{(1-a+a^2/2+b+b^2/2)^3}, 
\]
and this is not negated by the transformation.

However, we can still give a heuristic explanation for why Conjecture \ref{conj:odd} is true. Consider the functional equation that we know the invariant density must obey:
\begin{equation}\label{eq:hfunctional}
h(z)= \sum_{(\alpha+z)^{-1}\in T^{-1}z} \frac{1}{|\alpha +z|^4}h\left(\dfrac{1}{\alpha+z}\right), \ z\in K.
\end{equation}
If $z$ is treated as $x+y \ii$, then by taking appropriate real derivatives with respect to $x$ and $y$, we can use this to derive functional equations that should be satisfied by the various derivatives of $h(z)$. However, the more interesting aspect is the sum. The sum is over $(\alpha+z)^{-1}\in T^{-1}z$. All of the $\alpha$'s must belong to $\mathbb{Z}[\ii]$. Thus, if we let $z=-.5-.5\ii$, then $z+\mathbb{Z}[\ii]$ is a lattice of points that is symmetric with respect to $180$ degree rotation around the origin. As a consequence the possible values for $(\alpha+z)^{-1}\in K$ will also be symmetric with respect to $180$ degree rotation around the origin. 

In particular, the right-hand side of \eqref{eq:hfunctional} can be decomposed into pairs like
\[
\frac{1}{|\alpha+z|^4}h\left( \dfrac{1}{\alpha+z}\right) + \frac{1}{|-(\alpha+z)|^4}h\left( -\dfrac{1}{\alpha+z}\right) .
\]
We know that $h(z)=h(-z)$, because the invariant measure satisfies $\mu(A)=\mu(\mathbbm{i}A)$ and thus $\mu(A)=\mu(-A)$ by Theorem \ref{thm:Hensley}, so we expect odd derivatives of the above pair-sum to disappear. This holds perfectly well on all points $(a+z)^{-1}$ which are interior to one of the regions $K_{k,\ell}$, since here $h$ is analytic; however, it is not clear how to make this process work for those $(\alpha+z)^{-1}$ which lie on the boundary between regions, since $h$ is no longer analytic at these points. Moreover as $z=-.5-.5\ii$ lies on the boundary of $K$, such $(\alpha+z)^{-1}$ lying on the boundary between regions will exist. Thus, we can only describe this idea heuristically.

This heuristic cannot tell us the whole story. $\mathbb{Z}[\ii]-.5-.5\ii$ is not the only lattice that is symmetric with respect to $180$ degree rotation around the origin. $\mathbb{Z}[\ii]$ itself is symmetric in this way. So we would expect to see a similar cancellation in the odd derivatives of $h_{m,n}$ when starting our calculations from $z=0$ instead of $z=-.5-.5\ii$. 

However, in Table \ref{tab:V21table}, we calculated $h_{m,n}^{2,1}(0,0)$, since $0$ is on the boundary of $K_{2,1}$, and yet the odd terms are not clearly tending towards zero. We also calculated $h_{m,n}^{2,1}(-.5,-.5)$  in Table \ref{tab:V21table2}, but again the odd terms are not clearly tending towards zero. We note that $-.5-.5\ii$ is not in or on the border of $K_{2,1}$, but due to the similarity of $V_{2,1}$ and $V_{1,1}$, we thought it would be an interesting point to test regardless.

However, when we examined $V_{3,1}$, we again noticed that the odd $x$-derivatives tended to zero when we centered at $z=-.5$, see Table \ref{tab:V31table}.  As noted above, the natural symmetries of $V_{3,1}$ around the $y$-axis tell us that the $h$ should be even around this point in the $y$ direction, but our calculations also lead us to the following conjecture.

\begin{conj}\label{conj:odd}
We have $h_{m,n}^{3,1}(-.5,0)=0$ whenever $m$ or $n$ is odd.
\end{conj}
 
Heuristically we expect this to hold for the same reason as the previous conjecture.

\section{Other methods of calculating the invariant measure}\label{sec:othermethods}

In this section, we want to briefly discuss various different methods of calculating the density of the invariant measure and the benefits and problems with each. 

\subsection{The finite difference method}

It is well known that derivatives can be approximated by finite differences. So, for example, for a real differentiable function $f$ and sufficiently small $\epsilon$, we have
\begin{align*}
f'(x) &\approx \frac{f(x+\epsilon)-f(x)}{\epsilon}\\
f''(x) &\approx \frac{f(x+2\epsilon)-2f(x+\epsilon)+f(x)}{\epsilon^2}\\
f'''(x)&\approx \frac{f(x+3\epsilon)-3f(x+2\epsilon)+3f(x+\epsilon)-f(x)}{\epsilon^3}
\end{align*}
and so on.

To implement this, let us consider the problem of calculating, say $\frac{\partial^2}{\partial x^2} h(x+y \ii)$ at a point $x,y$. We pick a very small $\epsilon>0$ and let $U_1,U_2,U_3$ be balls of radius $\epsilon/2$ such that the center of $U_n$ is $x+y \ii+ (n-1)\epsilon$ for $n=1,2,3$. These sets will  be disjoint. Let $z_0$ be some point whose orbit under the Gauss map we expect to equidistribute over $K$, and let $Q$ be some extremely large integer.  

Now for $n=1,2,3$, let $u_n$ be the number of times $T^i z_0\in U_n$ as $i$ ranges from $0$ to $Q$. Therefore, $h(x+y \ii +(n-1)\epsilon)$ is approximated by the quantity $u_n/Q$. So our desired second derivative is approximated by
\[
\frac{u_3-2u_2+u_1}{Q\epsilon^2}.
\]

The main benefit of this method is in its simplicity. This could be implemented in only a few lines of code with very few variables being used.

The trade-off, naturally, is in accuracy. In order to have a good approximation for our coefficient, $\epsilon$ needs to be very small. But this makes the sets $U_n$ very small, and so $Q$ needs to be very large in order for $u_n$ to accurately approximate the density $h(x+y \ii +(n-1) \epsilon)$. Furthermore, even if the orbit of $z_0$ is equidistributed, it could ``approach" equidistribution very slowly; the values $u_n$ might converge to $\mu(U_n)$, but might do so very slowly, with lots of variance. This would encourage us to take $Q$ even larger.

\subsection{Functional equation method}

Recall that the Perron-Frobenius operator gave rise to the functional equation
\[
h(z)= \sum_{(a+z)^{-1}\in T^{-1}z} \frac{1}{|a+z|^4}h\left(\dfrac{1}{a+z}\right), \ z\in K.
\]
Let us choose a point $z_{k,\ell}\in K_{k,\ell}$ in each of the twelve regions and consider expanding $h(z')$ when $z'\in K_{k,\ell}$ by its Taylor series around the point $z_{k,\ell}$. For example, if $z\in K_{k,\ell}$, then the left-hand side becomes
\[
h(z)= \sum_{m,n=0}^\infty \frac{h_{m,n}^{k,\ell}(z_{k,\ell})}{m!n!}\Re(z-z_{k,\ell})^m\Im(z-z_{k,\ell})^n.
\]
For a fixed value of $z$, we can consider the quantities that look like
\[
\frac{1}{m!n!}\Re(z-z_{k,\ell})^m\Im(z-z_{k,\ell})^n
\]
as fixed values, so that the functional equation becomes a functional equation in the variables $h_{m,n}(z_{k,\ell})$ as $m,n$ range from $0$ to $\infty$ and $1\le k \le 3$, $1\le \ell\le 4$. Note further that this functional equation is linear in the variables $h_{m,n}(z_{k,\ell})$.

We can truncate this functional equation so that it has a finite number of terms by removing all summands involving $h_{m,n}(z_{k,\ell})$ with $m,n>N$ for some integer $N$. Precisely, this will have $12(N+1)^2$ terms. This truncation turns our exact functional equation into an approximate functional equation, but we ignore the errors and presume they are equal. 

Thus, if we take at least $12(N+1)^2+1$ different values of $z$, we obtain a system of linear (approximate functional) equations in the variables $h_{m,n}(z_{k,\ell})$, which we can then attempt to solve using a least squares approximation.

The advantage of using this method is that it is completely deterministic. Our other methods so far have all relied on having a point $z_0$ whose orbit under the Gauss map had good properties. We have no easy way of choosing a point with the desired property, so instead we chose a point essentially at random, since with probability $1$, we will choose a point with the desired property.

The disadvantage is that this method showed very poor convergence as $m,n$ increased. This happened because the coefficients multiplied to the $h_{m,n}$ shrink so rapidly that very large changes in $h_{m,n}$ would only alter the contribution of the overall term by a very small amount.

\subsection{Outline method}
 
In the method we used to give the results above, we make use of a flood-fill algorithm to fill in any ``holes" in \textbf{PixelArray}. As such, we only needed an approximation to the boundary of $V_{1,1}$ rather than the entirety of $V_{1,1}$. We took advantage of this in our variant of step 2, where we in essence only looked at points that were relatively close to the boundary. We attempted at one point to circumvent this entirely and compute those points that lie on the boundary itself. This was what led us to look at the admissible sequences in Section \ref{ref:admissible}, but we could not implement this in a way that was faster than the algorithm we ultimately used.

\end{document}